\documentclass[11pt,a4paper]{amsart}
\usepackage[T1]{fontenc}
\usepackage[utf8]{inputenc}
\usepackage[english]{babel}
\usepackage{lmodern,amsmath,amssymb,mathtools}
\usepackage[margin=29mm]{geometry}
\usepackage[expansion=false]{microtype}
\usepackage{enumitem}
\usepackage{booktabs,array}
\usepackage[hidelinks]{hyperref}
\hypersetup{
 pdftitle={Fourier--Stieltjes palettes on compact abelian groups},
 pdfauthor={Przemyslaw Ohrysko},
 pdfsubject={Exact Fourier--Stieltjes palettes on compact abelian groups},
 pdfkeywords={Fourier--Stieltjes transform, exact range, singular measure, Rajchman measure, compact abelian group, Wiener--Pitt set, Rudin--Shapiro polynomial}}
\allowdisplaybreaks[2]
\numberwithin{equation}{section}
\newtheorem{theorem}{Theorem}[section]
\newtheorem{lemma}[theorem]{Lemma}
\newtheorem{proposition}[theorem]{Proposition}
\newtheorem{corollary}[theorem]{Corollary}
\theoremstyle{definition}
\newtheorem{definition}[theorem]{Definition}

\theoremstyle{remark}
\newtheorem{remark}[theorem]{Remark}
\newcommand{\T}{\mathbb T}
\newcommand{\Z}{\mathbb Z}
\newcommand{\N}{\mathbb N}
\newcommand{\C}{\mathbb C}
\newcommand{\E}{\mathbb E}
\newcommand{\wh}{\widehat}
\newcommand{\norm}[1]{\lVert#1\rVert}
\newcommand{\HH}{\mathcal H}
\newcommand{\supp}{\operatorname{supp}}
\newcommand{\eps}{\varepsilon}
\newcommand{\one}{\mathbf 1}
\newcommand{\Tor}{\operatorname{Tor}}
\DeclareMathOperator{\RePart}{Re}

\title[Fourier--Stieltjes palettes]{Fourier--Stieltjes palettes on compact abelian groups}
\author{Przemys\l aw Ohrysko}
\address{University of Warsaw, Faculty of Mathematics, Informatics and Mechanics,
Institute of Mathematics, Banacha 2, 02-097 Warsaw, Poland}
\email{P.Ohrysko@mimuw.edu.pl}
\urladdr{https://orcid.org/0000-0002-2763-7031}
\date{November 16, 2026}

\begin{document}

\begin{abstract}
Let $G$ be an infinite compact abelian group and let
$K=\{0,a_1,a_2,\ldots\}$, where
$a_1>a_2>\cdots>0$ and $a_n\to0$.
We construct a nonzero real, even, purely singular Rajchman measure whose
Fourier--Stieltjes range is contained in $K$. We then classify the compact
abelian groups on which every positive null palette is exactly realizable
in each of four classes: purely singular measures, atomless purely singular
measures, strongly continuous purely singular measures, and Rajchman purely
singular measures. Exact realization by a purely singular measure is always
possible; atomless exact realization fails precisely for groups
topologically isomorphic to $\T\times F$ with $F$ finite; and universal
exact Rajchman realization is equivalent to $\Tor(\widehat G)$ being
infinite, or equivalently to $G/G_0$ being infinite. On $\T$ we also give
palette-dependent sufficient conditions for exact Rajchman realization.
Finally, if $U\subset\C$ is open and $0\in\overline U$, we construct a
continuous Rajchman measure with nonzero singular part whose
Fourier--Stieltjes range is contained in $U$, assuming metrizability of $G$
when $0\notin U$. This gives a negative answer to an open-range
absolute-continuity question of Ohrysko and Wojciechowski.
\end{abstract}

\maketitle

\section{Definitions, background, and main results}\label{sec:intro}

\subsection{Fourier--Stieltjes palettes}

Throughout, unless explicitly stated otherwise, $G$ is an
\emph{infinite} compact Hausdorff abelian group, $m_G$ is its normalized
Haar measure, and $\Gamma=\wh G$ is its discrete dual. This convention
removes only a degenerate case. If $G$ is finite, then every measure is
absolutely continuous with respect to $m_G$ and every Fourier--Stieltjes
transform has finite range; hence there are no nonzero purely singular
measures and no positive null palette can be realized exactly. We shall
therefore not repeat the word ``infinite'' in subsequent statements.

We use additive notation for both groups and evaluate characters
multiplicatively. For $\mu\in M(G)$,
\[
 \wh\mu(\gamma)=\int_G\overline{\gamma(x)}\,d\mu(x),
 \qquad \gamma\in\Gamma.
\]
The measure algebra $M(G)$ is equipped with convolution and total
variation norm. Its identity is the point mass $\delta_0$. We identify
$L^1(G)$ with the closed ideal of measures absolutely continuous with
respect to $m_G$.

A \emph{positive null palette} is a compact set
\begin{equation}\label{eq:palette}
 K=\{0,a_1,a_2,\ldots\},\qquad
 a_1>a_2>\cdots>0,\qquad a_n\longrightarrow0.
\end{equation}
There are two distinct range problems. A measure has \emph{range
contained in $K$} if $\wh\mu(\Gamma)\subset K$. It \emph{realizes $K$
exactly} if
\begin{equation}\label{eq:exact-palette}
 \wh\mu(\Gamma)=K.
\end{equation}
The distinction is essential: the general construction below selects a
subsequence of the prescribed values, whereas exact realization requires
every $a_n$ to occur at least once and also requires the value $0$ to
occur.

We write $M_0(G)$ for the Rajchman measures, that is,
$\wh\mu\in c_0(\Gamma)$. A measure is called \emph{continuous} if it is
atomless. It is \emph{strongly continuous} if
\begin{equation}\label{eq:strong-continuity}
 |\mu|(x+H)=0
\end{equation}
for every $x\in G$ and every closed subgroup $H\le G$ of infinite
index. It is \emph{even} if it is invariant under $x\mapsto-x$.
Singularity always means $|\mu|\perp m_G$, not merely that the singular
part is nonzero. We denote by $G_0$ the connected component of the
identity and by $\Tor\Gamma$ the torsion subgroup of $\Gamma$.

Although the main results below concern compact groups, the spectral
background is most naturally stated for an arbitrary locally compact
abelian group $H$. Let $\widehat H$ denote its Pontryagin dual. The
convolution algebra $M(H)$ of bounded regular Borel measures is a unital
commutative Banach algebra with identity $\delta_0$. For $\nu\in M(H)$ we
write
\[
 \sigma_{M(H)}(\nu)
 =\{\lambda\in\C:\nu-\lambda\delta_0
     \text{ is not invertible in }M(H)\}
\]
for the spectrum of $\nu$. If $\Delta(M(H))$ is the maximal ideal space,
Gelfand theory gives
\[
 \sigma_{M(H)}(\nu)
 =\{\varphi(\nu):\varphi\in\Delta(M(H))\}.
\]
Every $\gamma\in\widehat H$ determines the multiplicative functional
$\varphi_\gamma(\nu)=\widehat\nu(\gamma)$, and hence
\begin{equation}\label{eq:fourier-spectrum-inclusion}
 \overline{\widehat\nu(\widehat H)}
 \subset \sigma_{M(H)}(\nu).
\end{equation}
A measure has a \emph{natural spectrum} when equality holds in
\eqref{eq:fourier-spectrum-inclusion}. Following Zafran
\cite{Zafran1973}, we call such a measure a measure with natural spectrum.

The \emph{Wiener--Pitt phenomenon} on $H$ is the possibility that the
inclusion in \eqref{eq:fourier-spectrum-inclusion} is strict. This is
equivalent to the failure of the naive Wiener inversion criterion. Indeed,
if
\[
 \lambda\in\sigma_{M(H)}(\nu)
 \setminus\overline{\widehat\nu(\widehat H)},
\]
then $\nu-\lambda\delta_0$ is not invertible, although
\[
 \inf_{\gamma\in\widehat H}
 \bigl|\widehat\nu(\gamma)-\lambda\bigr|>0.
\]
Conversely, if a measure has Fourier--Stieltjes transform bounded away
from zero but is not invertible, then $0$ is a spectral value which is
invisible on $\widehat H$.

There is a sharp group-theoretic dichotomy. If $H$ is discrete, then
$M(H)=\ell^1(H)$ and the maximal ideal space of $M(H)$ is canonically
$\widehat H$; consequently every measure has a natural spectrum and the
Wiener--Pitt phenomenon does not occur. If $H$ is nondiscrete, the abstract
Wiener--Pitt theorem, proved in full generality by Williamson
\cite{Williamson1959}, provides a measure $\nu\in M(H)$ such that
\[
 \inf_{\gamma\in\widehat H}|\widehat\nu(\gamma)|>0,
 \qquad
 \nu\notin M(H)^{-1}.
\]
Thus the phenomenon occurs on every nondiscrete locally compact abelian
group. It originates in the work of Wiener and Pitt
\cite{WienerPitt1938}. Sreider's analysis of maximal ideals
\cite{Sreider1950}, the construction of additional multiplicative
functionals by Hewitt and Kakutani \cite{HewittKakutani1960}, and Graham's
Riesz-product proof \cite{Graham1974} are among the classical developments;
see also \cite{GrahamMcGehee1979,Rudin1962}.

The mechanism is encoded in the complicated Gelfand space of the measure
algebra. For every $\varphi\in\Delta(M(H))$, the restriction of $\varphi$
to the closed ideal $L^1(H)$ is either zero or evaluation at a unique
$\gamma\in\widehat H$. In the latter case $\varphi=\varphi_\gamma$ on all
of $M(H)$. Hence
\[
 \Delta(M(H))\setminus h(L^1(H))
 =\{\varphi_\gamma:\gamma\in\widehat H\},
 \qquad
 h(L^1(H))
 =\{\varphi:\varphi|_{L^1(H)}=0\}.
\]
The first part is the visible copy of the dual group, whereas the hull
$h(L^1(H))$ consists of the extra multiplicative functionals responsible
for spectral values not detected by the Fourier--Stieltjes transform. For
discrete $H$ the hull is empty; for nondiscrete $H$ it is precisely this
additional part of $\Delta(M(H))$ that makes the Wiener--Pitt phenomenon
possible. The size of this additional part is not merely a heuristic indication of
complexity. Ohrysko, Wojciechowski and Graham proved that, for every
nondiscrete locally compact abelian group $H$, both $\Delta(M_0(H))$ and
$\Delta(M(H))$ contain continuum many pairwise disjoint nonempty open sets
and are therefore nonseparable
\cite[Theorem~11 and Corollary~12]{OWG2016}. They also proved that the
spectrum in $M(H)$, and already in $M_0(H)$, is not countably generated
\cite[Fact~17]{OWG2016}: for neither algebra $A$ is there a countable set
$E\subset\Delta(A)$ such that
\[
 \sigma_A(a)=\overline{\{\varphi(a):\varphi\in E\}}
 \qquad(a\in A).
\]
This gives a precise topological manifestation of the abundance of spectral
values not detected on $\widehat H$. Returning to the compact group $G$
fixed above, we use the same notation $\sigma_{M(G)}$ and the same notion
of natural spectrum.

A compact set $E\subset\C$ is called a Wiener--Pitt set if every measure
whose Fourier--Stieltjes range is contained in $E$ has a natural spectrum.
For each prescribed norm bound, the results of
\cite[Theorems~2 and~12]{OSW2026} show that sufficiently rapidly decreasing
positive palettes force the convolution square of every strongly
continuous measure in that norm ball to be absolutely continuous; such
measures also have natural spectra. This is a norm-bounded version of the
Wiener--Pitt property for a specified class of measures, rather than the
unrestricted property in the preceding definition. Exact realization asks
a different question: whether all values of the palette can occur
simultaneously for a singular measure.

Classical range-rigidity results and the related almost-idempotent theory
are recalled in the historical context below. The present problem is
complementary: zero is the only accumulation point, but every nonzero value
is prescribed in advance and, in the exact problem, must actually occur.

\subsection{Historical context}

The problem belongs to several classical lines of commutative harmonic
analysis. The first concerns singular measures with small transforms.
For each fixed $\eta>0$, Wiener and Wintner
\cite{WienerWintner1938} constructed a singular probability measure
$\mu_\eta$ on the circle with
\[
 \wh{\mu_\eta}(n)=O(|n|^{-1/2+\eta}).
\]
For every strictly increasing function
$\rho:\N\to(0,\infty)$ with $\rho(n)\to\infty$, Schaeffer
\cite{Schaeffer1939} constructed a singular probability measure $\mu_\rho$
with
\[
 \wh{\mu_\rho}(n)=O\!\left(\frac{\rho(|n|)}{\sqrt{|n|}}\right).
\]
Hewitt and Ritter
\cite{HewittRitter1976} obtained refined near-$L^2$ integrability
results on locally compact abelian groups; see also
\cite[Theorems~3.1--3.2]{HewittRitter1977}. Brown and Hewitt
\cite{BrownHewitt1980} developed further constructions of continuous
singular measures with small Fourier--Stieltjes transforms. Kahane
\cite{Kahane2007} used products of Rudin--Shapiro blocks to obtain
singular measures satisfying very strong Orlicz summability conditions.
Those results control the decay or size distribution of the transform;
the present problem prescribes its countable range, including the
accumulation point zero.

A second line concerns regular convolution powers. Hewitt and Zuckerman
\cite{HewittZuckerman1966,HewittZuckerman1967} constructed continuous
singular measures whose convolution squares have densities in
$L^p(\T)$ for every $1\le p<\infty$. Stromberg obtained large families of such measures
\cite{Stromberg1968}. Saeki developed the theory for general groups and
support restrictions \cite{Saeki1977}; on the circle he obtained a
singular probability measure whose convolution square has a uniformly
convergent Fourier series \cite{Saeki1980}. Riesz-product proofs and
extensions were given by Karanikas and Koumandos
\cite{KaranikasKoumandos1991} and Bisbas \cite{Bisbas1999}.
Andersson \cite{Andersson2002} analyzed three known constructions with
uniformly bounded $L^1$ Fourier partial sums and showed that their
convolution squares have $L^2$ densities. K\"orner \cite{Korner2008} combined
regularity of the square with control of the Hausdorff dimension of the
support. A later real-even flat-density assertion is discussed critically
in Remark~\ref{rem:el-abdalaoui}; it is not used in any of our arguments.

A third line studies the structural information encoded by the range.
In the case of an isolated zero, Glicksberg's theorem
\cite[Theorem~1.1, pp.~59--60]{Glicksberg1972} identifies the structure of
the corresponding part of the measure, while Glicksberg and Wik
\cite{GlicksbergWik1972} obtained range-inclusion results for suitable
parts of measures.
The almost-idempotent theory of de Leeuw--Katznelson
\cite{DeLeeuwKatznelson1970} and Pigno--Smith
\cite{PignoSmith1980} quantifies related rigidity. Sarnak
\cite{Sarnak1980} obtained spectral conclusions for convolution
operators on $L^p(\T)$ from restrictions on transform values; those
operator spectra must be distinguished from the spectrum in $M(\T)$.
The Wiener--Pitt-set program in \cite{OSW2026,OW2014} asks when the value
set alone forces a natural spectrum or regular convolution powers.

Finally, frequency-domain interpolation is different from value-domain
prescription. Wells \cite{Wells1973} characterized Sidon sets through
restrictions of transforms of continuous measures, and Wojciechowski
\cite{Wojciechowski2018} studied zero sets of Fourier transforms of
singular measures. In those problems one prescribes behavior on a set of frequencies. Here
we constrain the range of the entire transform rather than assign values
at preselected frequencies. Standard background can be found in
\cite{GrahamMcGehee1979,Katznelson2004,Rudin1962}.

The combinatorial core of our proof is the Rudin--Shapiro recursion
\cite{Rudin1959,Shapiro1951}. Its quadratic path-sign description is
closely related to the algebraic description of Golay complementary
sequences in \cite[Theorem~3]{DavisJedwab1999}. The Boolean selector
built from a paired colouring of the path is what turns those signs into exactly
prescribed coefficient values.

\subsection{The general contained-range theorem}

The first result is the general construction proved in
Section~\ref{sec:general-construction}.

\begin{theorem}\label{thm:main}
Let $G$ be a compact Hausdorff abelian group and let $K$ be as
in~\eqref{eq:palette}. For every $\eps>0$ there is a real even measure
$\mu\in M(G)$ such that
\begin{gather}
 \mu\perp m_G,\qquad \wh\mu(0)=a_1,\qquad
 \wh\mu(\Gamma)\subset K,\label{eq:main-range}\\
 a_1\le\norm\mu\le a_1+\eps,\qquad
 \wh\mu\in\bigcap_{r>2}\ell^r(\Gamma),\qquad
 \mu*\mu\in\bigcap_{1\le p<\infty}L^p(G).
 \label{eq:main-regularity}
\end{gather}
Moreover, $\mu^{*3}$ has a continuous density with an absolutely
convergent Fourier series. The measure belongs to $M_0(G)$, is strongly
continuous, and has a natural spectrum with spectral radius $a_1$.
\end{theorem}

No condition is imposed on the ratios $a_{n+1}/a_n$, and no
metrizability assumption is made. The endpoint $\ell^2$ is impossible:
a finite measure with square-summable Fourier coefficients is represented
by an $L^2$ density. The construction combines Rudin--Shapiro blocks,
Boolean selectors on a coloured path, finite Fourier averaging, and exact
recovery of the approximating polynomials.

\begin{corollary}\label{cor:norm-cap}
Let $G$ be a compact Hausdorff abelian group and let $K$ be as
in~\eqref{eq:palette}. For every $R>0$ there are an index $n_0$ and a
nonzero real even purely singular strongly continuous measure $\mu\in M_0(G)$ with
$\norm\mu<R$ and
\[
 \wh\mu(\Gamma)\subset\{0,a_{n_0},a_{n_0+1},\ldots\}\subset K,
 \qquad \wh\mu(0)=a_{n_0}.
\]
Its transform belongs to every $\ell^r(\Gamma)$, $r>2$, its convolution
square belongs to $L^p(G)$ for every $1\le p<\infty$, and its third convolution
power has a continuous density with an absolutely convergent Fourier
series. It has a natural spectrum with spectral radius $a_{n_0}$.
\end{corollary}

\begin{proof}
Choose $n_0$ with $a_{n_0}<R$ and apply Theorem~\ref{thm:main} to
$\{0,a_{n_0},a_{n_0+1},\ldots\}$ with
$0<\eps<R-a_{n_0}$. The spectral radius is the first value of this
chosen tail, rather than the first value of the original palette.
\end{proof}

\subsection{Exact palettes and the structure of the group}

We use four universal properties.

\begin{definition}\label{def:universal-properties}
For a compact abelian group $G$ let:
\begin{enumerate}[label=\textup{(\roman*)}]
\item $\mathcal U_{\mathrm s}(G)$ mean that every positive null palette
is exactly realized by a purely singular measure;
\item $\mathcal U_{\mathrm c}(G)$ mean that the realizing measure may
also be required to be atomless;
\item $\mathcal U_{\mathrm{sc}}(G)$ mean that it may be required to be
strongly continuous;
\item $\mathcal U_0(G)$ mean that it may be required to belong to
$M_0(G)$.
\end{enumerate}
In all four cases the measure is required to be purely singular. It may
be complex unless the theorem explicitly gives a real even realization.
\end{definition}

The four properties satisfy
\[
 \mathcal U_0(G)\Longrightarrow\mathcal U_{\mathrm{sc}}(G)
 \Longrightarrow\mathcal U_{\mathrm c}(G)
 \Longrightarrow\mathcal U_{\mathrm s}(G).
\]
The following classification is proved in
Section~\ref{sec:exact-palettes}.

\begin{theorem}[Classification of universal exact realization]
\label{thm:classification-overview}
Let $G$ be a compact abelian group.
\begin{enumerate}[label=\textup{(\alph*)}]
\item $\mathcal U_{\mathrm s}(G)$ always holds.
\item $\mathcal U_{\mathrm c}(G)$ holds if and only if $G$ is not
topologically isomorphic to $\T\times F$ for a finite group $F$.
\item The following are equivalent:
\[
 \mathcal U_{\mathrm{sc}}(G),\qquad
 \mathcal U_0(G),\qquad
 |\Tor\wh G|=\infty,\qquad
 |G/G_0|=\infty.
\]
In the affirmative case of part~\textup{(c)}, for every $\eps>0$ the
realizing measure can be chosen real and even with mass $a_1$ and
$a_1\le\norm\mu\le a_1+\eps$.
\end{enumerate}
\end{theorem}

The three boundaries arise from different mechanisms. Infinite torsion
in the dual permits forcing an arbitrarily long finite tranche of the
palette at a cost depending only on its largest value. An infinite-order
character yields a completely explicit telescoping construction. The
negative statements follow from quantitative gap theorems for Fourier--
Stieltjes transforms of continuous or strongly continuous measures.

The classification concerns universal realization of \emph{every}
palette. On the circle, where atomless universal realization fails, many
individual palettes can nevertheless be exactly realized by singular
Rajchman measures.

\begin{theorem}[Exact palettes on the circle]\label{thm:circle-overview}
Let $K$ be as in~\eqref{eq:palette}.
\begin{enumerate}[label=\textup{(\alph*)}]
\item If
\[
 \sum_{n\ge1}\frac{a_n}{n}<\infty,
\]
then $K$ is exactly realized by a real even purely singular Rajchman
measure on $\T$.
\item The same conclusion holds if $(a_n)$ is discretely convex:
$a_n-2a_{n+1}+a_{n+2}\ge0$ for all $n$. No prescribed decay rate is
needed in this case.
\end{enumerate}
\end{theorem}

\subsection{Applications and questions from the literature}

The method also has consequences which do not follow formally from the
exact-realization results. First, the open-range question posed in
\cite[Section~7, item~2]{OW2014} has a negative answer.

\begin{theorem}\label{thm:complex}
Let $G$ be a compact abelian group and let $(z_n)$ be a sequence of
nonzero complex numbers tending to zero. For every $R>0$ there is a strongly
continuous $\nu\in M_0(G)\setminus L^1(G)$ such that
\[
 \norm\nu<R,\qquad
 \wh\nu(\Gamma)\subset\{0,z_1,z_2,\ldots\},\qquad
 \wh\nu\in\bigcap_{r>2}\ell^r(\Gamma).
\]
The second and third convolution powers have the same regularity as in
Theorem~\ref{thm:main}, and $\nu$ has a natural spectrum.
\end{theorem}

\begin{theorem}\label{thm:open}
Let $G$ be a compact abelian group and let $U\subset\C$ be open
with $0\in\overline U$. If $0\notin U$, assume that $G$ is metrizable.
For every $R>0$ there is a strongly continuous
$\nu\in M_0(G)\setminus L^1(G)$ with
\[
 \norm\nu<R,\qquad \wh\nu(\Gamma)\subset U,
 \qquad \wh\nu\in\bigcap_{r>2}\ell^r(\Gamma),
\]
and with the same convolution-power and natural-spectrum conclusions.
If $0\in U$, no metrizability assumption is needed.
\end{theorem}

Second, a phase-valued selector gives a reflection-invariant complex
analogue of Saeki's construction. On $\T$, uniform convergence is understood
with respect to the ordinary symmetric Fourier partial sums. This must be
distinguished from the unconditional finite-subset norm displayed in
\cite{GuptaHare}; under that stronger interpretation no nonzero singular
measure is possible.

\begin{theorem}\label{thm:uniform-square}
Let $A>0$, $\eps>0$, and let $(a_n)$ be a positive null sequence. There
are a subsequence $c_j=a_{n_j}$ and a nonzero even complex measure
$\nu\in M_0(\T)$ such that
\begin{gather}
 \nu\perp m_\T,\qquad \wh\nu(0)=A,
 \qquad A\le\norm\nu\le A+\eps,\label{eq:uniform-measure}\\
 \wh\nu(\Z)\subset\{0,A\}\cup
 \bigcup_{j\ge1}c_j\{1,-1,i,-i\},
 \qquad \wh\nu\in\bigcap_{r>2}\ell^r(\Z).
 \label{eq:uniform-range}
\end{gather}
The Fourier series of $\nu*\nu$ converges uniformly but not absolutely. In particular,
$\nu$ is strongly continuous and has a natural spectrum.
\end{theorem}

This answers the circle question for even complex measures when uniform
convergence is understood in terms of the ordinary symmetric Fourier partial
sums. Complex phases are essential under the evenness requirement.
Proposition~\ref{prop:real-even-obstruction} shows more generally
that a real even measure whose convolution square has an essentially
bounded density must itself have an $L^2$ density. Thus a nonzero purely
singular real even measure cannot have a bounded square density, let alone
a uniformly convergent convolution-square Fourier series. The apparently
conflicting real-even flat-density assertion in \cite{ElAbdalaoui2023} is
discussed in Remark~\ref{rem:el-abdalaoui} and is not used as a precedent.

\subsection{Organization and provenance of the arguments}

Section~\ref{sec:general-construction} contains the complete proof of
Theorem~\ref{thm:main}. Section~\ref{sec:exact-palettes} develops exact
realization and proves Theorem~\ref{thm:classification-overview}, as
well as the circle criteria in Theorem~\ref{thm:circle-overview}.
Section~\ref{sec:applications} proves Theorems~\ref{thm:complex},
\ref{thm:open}, and~\ref{thm:uniform-square}, and relates the results to
questions on Wiener--Pitt sets and regular convolution powers.

The block construction is presented self-containedly in
Section~\ref{sec:general-construction}. The exact-palette classifications
combine that construction with the classical gap theorem in
\cite[Theorem~1.4.1]{GrahamMcGehee1979}.
\section{The general construction on compact abelian groups}\label{sec:general-construction}

\subsection{Quotients and the three model groups}\label{sec:reduction}

We begin with the reduction which permits a sequential construction
even when the original group is not metrizable.

\begin{lemma}[Haar lifting]\label{lem:lift}
Let $\pi:G\to H$ be a continuous surjective homomorphism of compact
abelian groups. There is a positive linear map $L_\pi:M(H)\to M(G)$
such that
\[
 \pi_*L_\pi\nu=\nu,\qquad
 \norm{L_\pi\nu}=\norm\nu,\qquad
 |L_\pi\nu|=L_\pi|\nu|.
\]
The Fourier transform of $L_\pi\nu$ agrees with $\wh\nu$ on
$\pi^*\wh H\subset\wh G$ and vanishes outside this subgroup.
The map preserves singularity, reality, evenness, the Rajchman
condition, and every $\ell^r$ norm of the transform. Furthermore,
\[
 L_\pi(\nu*\rho)=L_\pi\nu*L_\pi\rho,
 \qquad L_\pi(fm_H)=(f\circ\pi)m_G.
\]
In the second identity, the pullback $f\mapsto f\circ\pi$ preserves the
$L^p$ norm for $1\le p\le\infty$.
\end{lemma}

\begin{proof}
Put $N=\ker\pi$. For $f\in C(G)$ define $Af\in C(H)$ by
\[
 Af(\pi(x))=\int_N f(x+y)\,dm_N(y).
\]
The right-hand side is continuous and constant on the fibres of
$\pi$, so the definition is valid. Define the measure $L_\pi\nu$ by
\[
 \int_G f\,dL_\pi\nu=\int_H Af\,d\nu.
\]
Positivity and $\norm A=1$ give
$\norm{L_\pi\nu}\le\norm\nu$ and
$|L_\pi\nu|\le L_\pi|\nu|$.
For $g\in C(H)$ we have $A(g\circ\pi)=g$, whence
$\pi_*L_\pi\nu=\nu$. This gives the reverse norm inequality.
The two positive measures in the variation inequality have the same
total mass, proving equality.

Averaging a character over $N$ gives zero unless it is trivial on $N$;
in the latter case the character factors through $H$.
This proves the Fourier formula. Fourier uniqueness then proves the
convolution identity. The formula for densities follows from Haar
integration over the fibres, first for continuous functions and then
by $L^1$ approximation. Since $\pi_*m_G=m_H$, pullback is isometric
on every $L^p$.

If $|\nu|$ is concentrated on a Borel set $E$ with $m_H(E)=0$,
the variation identity and the pushforward identity show that
$|L_\pi\nu|$ is concentrated on $\pi^{-1}(E)$.
This set has Haar measure zero. The remaining assertions follow
from the Fourier formula or directly from the definition.
\end{proof}

\begin{lemma}\label{lem:algebra}
Every infinite discrete abelian group contains a subgroup of one of
the following types:
\begin{enumerate}[label=\textnormal{(\roman*)},leftmargin=2.2em]
\item $\Z$;
\item a countable torsion group of unbounded exponent;
\item $\bigoplus_{n\ge1}\Z/p\Z$ for some prime $p$.
\end{enumerate}
\end{lemma}

\begin{proof}
An element of infinite order gives the first case. In a torsion group
of unbounded exponent, choose elements with orders tending to infinity
and take the subgroup they generate. It is countable and gives the
second case.

Otherwise the group has finite exponent. Its decomposition into
finitely many primary components has an infinite component $A$
annihilated by $p^s$ for some prime $p$ and some $s\ge1$.
If $A[p]=\{x\in A:px=0\}$ were finite, every fibre of multiplication
by $p$ would have at most $|A[p]|$ elements. Starting from $p^sA=0$
and proceeding backwards gives $|A|\le|A[p]|^s$, a contradiction.
Thus $A[p]$ is an infinite vector space over $\mathbb F_p$, and
contains a countably infinite linearly independent set.
\end{proof}

By Pontryagin duality, a subgroup $\Lambda\subset\Gamma$ gives a
quotient map $G\to\wh\Lambda$. The preceding lemmas reduce the
construction to three metrizable models:
\begin{equation}\label{eq:models}
 \T,\qquad
 \wh\Lambda\ \text{with $\Lambda$ countable torsion of unbounded exponent},
 \qquad (\Z/p\Z)^{\N}.
\end{equation}
The last two models have countable torsion duals. Every finitely
generated subgroup of a torsion abelian group is finite.

\subsection{Paired path colours and Boolean selectors}\label{sec:selectors}

The complementary-polynomial recursion used below is classical
\cite{Rudin1959,Shapiro1951}. For a related product construction using
Rudin--Shapiro blocks in the Walsh setting, see \cite{Kahane2007}; compare
also Section~6 of \cite{OW2014}. Here we combine it
with Boolean selectors so that every nonzero coefficient of each new
block equals one prescribed positive value.

Let $\mathbb U$ be the circle group of complex numbers of modulus one,
or the group of $h$th roots of unity, $h\ge2$, with normalized Haar
measure. All expectations and norms on products of $\mathbb U$ in
this section refer to independent Haar variables.

\subsubsection{The full path and an isolated vertex}

For $m\ge1$ define
\begin{equation}\label{eq:RS}
 R_m(z_1,\ldots,z_m)
 =\sum_{x\in\{0,1\}^m}
 (-1)^{\sum_{j=1}^{m-1}x_jx_{j+1}}z_1^{x_1}\cdots z_m^{x_m}.
\end{equation}
Let $T_m$ be the same polynomial with an additional coefficient
factor $(-1)^{x_m}$, and put $R_0=T_0=1$.
Splitting the sum according to the last bit gives
\[
 R_m=R_{m-1}+z_mT_{m-1},\qquad
 T_m=R_{m-1}-z_mT_{m-1}.
\]
Consequently
\[
 |R_m|^2+|T_m|^2=2^{m+1}
\]
pointwise. The monomials $z^x$, $x\in\{0,1\}^m$, are distinct
characters even when $h=2$. Orthogonality therefore gives
\begin{equation}\label{eq:RS-bounds}
 Z_m=2^{-m/2}R_m,\qquad
 \norm{Z_m}_2=1,\qquad \norm{Z_m}_\infty\le\sqrt2.
\end{equation}
This is the complementary-polynomial identity in the form needed here.

An $L^1/L^2$ ratio uniformly below one will be required only for a
single isolated vertex. Put $V=|1+z|/\sqrt2$, with $z$ Haar distributed
on $\mathbb U$. Direct expansion, using the first two moments of $z$,
gives
\[
 0\le V\le\sqrt2,\qquad \E V^2=1,
 \qquad
 \E V^4=
 \begin{cases}
 2,&h=2,\\
 3/2,&h\ge3\text{ or }\mathbb U=\T.
 \end{cases}
\]
It follows that
\[
 \frac12\le\E(V^2-1)^2
 \le(1+\sqrt2)^2\E(V-1)^2
 =2(1+\sqrt2)^2(1-\E V).
\]
Hence, uniformly in $\mathbb U$,
\begin{equation}\label{eq:q}
 \E\frac{|1+z|}{\sqrt2}\le q,
 \qquad q:=1-\frac1{4(1+\sqrt2)^2}<1.
\end{equation}

\subsubsection{Proper sets of colours}

Fix $d\ge1$, write $[d]=\{1,\ldots,d\}$, and colour the consecutive
edges of the path on $m$ vertices with the repeating pattern
\begin{equation}\label{eq:colours}
 1,1,2,2,\ldots,d,d,1,1,2,2,\ldots,d,d,\ldots.
\end{equation}
If $c(j)$ is the colour of the edge joining vertices $j$ and $j+1$,
put
\[
 \sigma_i(x)=(-1)^{\sum_{j:c(j)=i}x_jx_{j+1}}
 \quad(1\le i\le d),
 \qquad
 R_{m,S}(z)=\sum_{x\in\{0,1\}^m}
       \prod_{i\in S}\sigma_i(x)z^x
\]
for $S\subset[d]$. In particular, $R_{m,[d]}=R_m$.

\begin{lemma}\label{lem:proper}
For every proper subset $S\subsetneq[d]$,
\begin{equation}\label{eq:proper}
 2^{-m/2}\norm{R_{m,S}}_1\le q^{r(m,d)},
 \qquad r(m,d)=\left\lfloor\frac{m-1}{2d}\right\rfloor.
\end{equation}
\end{lemma}

\begin{proof}
Delete all edges whose colours are outside $S$, and fix one missing
colour. Each complete period of $2d$ edges contains two consecutive
deleted edges of that colour. Their common vertex is isolated.
These vertices are distinct in different periods, so there are at
least $r(m,d)$ isolated vertices.

The remaining graph is a disjoint union of paths, including singleton
vertices. The polynomial $R_{m,S}$ factors as a product of their path
polynomials in disjoint sets of variables. After normalization by
$2^{-v/2}$, where $v$ is the number of vertices of a component, every
factor has $L^2$ norm one, hence $L^1$ norm at most one.
Each distinguished singleton factor has $L^1$ norm at most $q$ by
\eqref{eq:q}. Independence of the variable sets gives the claim.
\end{proof}

The repetition of colours in pairs is useful for groups of exponent two.
For example, $1+z_1+z_2-z_1z_2$ has constant modulus $2$ when
$z_1,z_2\in\{-1,1\}$. Thus an $L^1/L^2$ ratio uniformly below one cannot be
assumed for every nontrivial path. Isolated vertices provide the
uniform bound needed in~\eqref{eq:proper}.

For a function $F:\{-1,1\}^d\to\C$ and $S\subset[d]$, we use the
normalized Walsh coefficient
\begin{equation}\label{eq:Walsh-normalization}
 \wh F(S)
 =2^{-d}\sum_{s\in\{-1,1\}^d}F(s)\prod_{i\in S}s_i.
\end{equation}
With this normalization, Parseval's identity reads
\begin{equation}\label{eq:Walsh-Parseval}
 \sum_{S\subset[d]}|\wh F(S)|^2
 =2^{-d}\sum_{s\in\{-1,1\}^d}|F(s)|^2.
\end{equation}

\subsubsection{Selectors with coefficients zero and one}

\begin{lemma}\label{lem:selector}
For every $b\in2^{-d}\Z\cap[-1/2,1/2]$, there is a polynomial
\[
 B_b(z)=\sum_{x\in\{0,1\}^m}c_xz^x,\qquad c_x\in\{0,1\},
\]
such that
\begin{equation}\label{eq:selector}
 2^{-m/2}\norm{B_b-bR_m}_1
 \le2^{d/2}q^{r(m,d)}.
\end{equation}
The selector for $b=0$ may be chosen to be zero.
\end{lemma}

\begin{proof}
On the Boolean cube $\{-1,1\}^d$, choose a function
$F_b:\{-1,1\}^d\to\{0,1\}$ whose full Walsh coefficient is $b$.
To see that this is possible, the two parity classes each have
$2^{d-1}$ points, and
\[
 \wh F_b([d])=2^{-d}
 \bigl(\#\{\text{chosen points of positive parity}\}
       -\#\{\text{chosen points of negative parity}\}\bigr).
\]
Every permitted value of $b$ is obtained by selecting the appropriate
number of points in one class and none in the other.
For $b=0$ choose $F_b=0$.

Set
\[
 B_b(z)=\sum_x F_b(\sigma_1(x),\ldots,\sigma_d(x))z^x.
\]
Expanding $F_b$ in Walsh characters yields the exact identity
\[
 B_b=\sum_{S\subset[d]}\wh F_b(S)R_{m,S}.
\]
No distributional assertion about the vector
$(\sigma_1(x),\ldots,\sigma_d(x))$ is needed: this is a pointwise
algebraic identity. Parseval and Cauchy--Schwarz give
\[
 \sum_{S\subset[d]}|\wh F_b(S)|
 \le2^{d/2}\Bigl(\sum_S|\wh F_b(S)|^2\Bigr)^{1/2}
 \le2^{d/2}.
\]
Separate the full Walsh term $bR_m$ and apply Lemma~\ref{lem:proper}
to the remaining terms.
\end{proof}

\subsection{A flat block on a product group}\label{sec:product}

For a measurable function $f$ on a probability space write
\[
 \HH(f)=\int\sqrt{|f|}.
\]
Let $Y$ be a compact abelian group, let $A>0$, and let
\begin{equation}\label{eq:old-P}
 P(y)=\sum_{\lambda\in E}p_\lambda\lambda(y),\qquad
 E=-E,\quad D=|E|,\quad
 0\le p_\lambda\le A,\quad p_{-\lambda}=p_\lambda.
\end{equation}
The finite set $E$ may include characters with zero coefficients.
The polynomial $P$ is real and even.
Introduce independent Haar variables $z_1,\ldots,z_m,w\in\mathbb U$.
If $\mathbb U=\{-1,1\}$, assume also that all characters in $E$
are real-valued. Put
\[
 \kappa=
 \begin{cases}
 1,&\mathbb U=\{-1,1\},\\
 2,&\text{otherwise}.
 \end{cases}
\]

\begin{lemma}\label{lem:product}
There is an absolute constant
\[
 \theta=1-\frac1{1024}<1
\]
with the following property. For every $P$ as in~\eqref{eq:old-P}
and every $\eta>0$, for \emph{every sufficiently small} $a>0$
one can choose $m\ge1$ and a real even trigonometric polynomial
$\mathcal Q$ on $Y\times\mathbb U^m\times\mathbb U$ such that:
\begin{enumerate}[label=\textnormal{(\roman*)},leftmargin=2.2em]
\item every nonzero Fourier coefficient of $\mathcal Q$ equals $a$;
\item its Fourier support is disjoint from the old characters
$(\lambda,0,\ldots,0)$;
\item
\begin{equation}\label{eq:product-estimates}
 \norm{P+\mathcal Q}_1\le\norm P_1+\eta,
 \qquad
 \HH(P+\mathcal Q)\le\theta\HH(P)+\eta;
\end{equation}
\item
\begin{equation}\label{eq:block-count}
 \#\supp\wh{\mathcal Q}\le\frac{DA^2}{4a^2}.
\end{equation}
\end{enumerate}
Here $P$ is viewed as a function on the product through the first
coordinate.
\end{lemma}

\begin{proof}
Choose dyadic coefficients
$b_\lambda\in2^{-d}\Z\cap[0,1/2]$ satisfying $b_{-\lambda}=b_\lambda$ and
\begin{equation}\label{eq:dyadic}
 A_d=\sum_{\lambda\in E}b_\lambda\lambda,
 \qquad
 \norm{A_d-P/(2A)}_\infty\le D2^{-d}.
\end{equation}
For example, round each $p_\lambda/(2A)$ down to this dyadic grid.
The symmetry relations $b_{-\lambda}=b_\lambda$ and $E=-E$ make $A_d$
real and even.
Use the same Boolean selector whenever the values $b_\lambda$ agree.
In particular, coefficients with $p_\lambda=0$ use the zero selector.

For $a>0$ small enough that the following integer is at least one, put
\begin{equation}\label{eq:m-choice}
 m=\left\lfloor2\log_2\frac{A}{2\kappa a}\right\rfloor,
 \qquad M=2^m,\qquad s=a\sqrt M,\qquad
 \beta=\frac{\kappa s}{2A}.
\end{equation}
Then
\begin{equation}\label{eq:beta}
 \frac{A}{2\kappa\sqrt2}<s\le\frac{A}{2\kappa},
 \qquad \frac1{4\sqrt2}<\beta\le\frac14.
\end{equation}
With $Z=M^{-1/2}R_m$, define
\[
 X=
 \begin{cases}
 wZ,&\kappa=1,\\
 \RePart(wZ),&\kappa=2.
 \end{cases}
\]
When $\kappa=1$, the polynomial $Z$ is real. In the other case,
$\E w=\E w^2=0$. Together with~\eqref{eq:RS-bounds} this gives
\begin{equation}\label{eq:X-moments}
 \E X=0,\qquad |X|\le\sqrt2,
 \qquad
 \E X^2=
 \begin{cases}1,&\kappa=1,\\1/2,&\kappa=2.\end{cases}
\end{equation}
In particular $|\beta X|<1/2$ and $1+\beta X>0$.
For $|u|\le1/2$ the elementary inequality
\[
 \sqrt{1+u}\le1+\frac u2-\frac{u^2}{16}
\]
follows, for example, from Taylor's formula and
$\frac{d^2}{du^2}\sqrt{1+u}=-1/[4(1+u)^{3/2}]$.
Consequently~\eqref{eq:beta} and~\eqref{eq:X-moments} imply
\begin{equation}\label{eq:contraction}
 \E\sqrt{1+\beta X}
 \le1-\frac{\beta^2\E X^2}{16}
 \le\theta.
\end{equation}
Independence from $y$ now yields the ideal estimates
\begin{equation}\label{eq:ideal}
 \norm{P(1+\beta X)}_1=\norm P_1,
 \qquad
 \HH(P(1+\beta X))\le\theta\HH(P).
\end{equation}
The equality and the inequality above use the positivity of $1+\beta X$,
not that of $P$.

Define
\[
 S(y,z)=\sum_{\lambda\in E}\lambda(y)B_{b_\lambda}(z),
 \qquad
 \mathcal Q(y,z,w)=
 \begin{cases}
 awS(y,z),&\kappa=1,\\
 a\bigl(wS(y,z)+\overline{wS(y,z)}\bigr),&\kappa=2.
 \end{cases}
\]
In the first case all characters are real, so $\mathcal Q$ is real
and even. In the second case conjugation supplies the negative
frequencies with the same coefficients. The two copies are disjoint,
because the exponents $1$ and $-1$ in the $w$ coordinate give distinct
characters. Thus every nonzero coefficient is exactly $a$, and the
nonzero $w$ exponent separates the support from every old character.

The full Walsh components sum to $\kappa sA_dX$.
The selector estimate, followed by~\eqref{eq:dyadic}, gives
\begin{align}
 \norm{\mathcal Q-\beta PX}_1
 &\le\kappa s\sqrt2\norm{A_d-P/(2A)}_\infty
       +\kappa sD2^{d/2}q^{r(m,d)}\notag\\
 &\le\frac{A}{\sqrt2}D2^{-d}
       +\frac A2D2^{d/2}q^{r(m,d)}=:\delta.\label{eq:product-error}
\end{align}
For integrable $f,g$ on a probability space,
\begin{equation}\label{eq:H-continuity}
 |\HH(f)-\HH(g)|
 \le\int\sqrt{|f-g|}\le\sqrt{\norm{f-g}_1}.
\end{equation}
Choose $d$ first so that the first term of~\eqref{eq:product-error}
is sufficiently small. With this $d$ fixed, $m\to\infty$ as
$a\downarrow0$, so the second term becomes sufficiently small for
\emph{every} sufficiently small $a$. We can therefore ensure
$\delta\le\min(\eta,\eta^2)$. Equations~\eqref{eq:ideal}
and~\eqref{eq:H-continuity} prove~\eqref{eq:product-estimates}.
Finally,
\[
 \#\supp\wh{\mathcal Q}\le\kappa D2^m
 \le\frac{DA^2}{4\kappa a^2}\le\frac{DA^2}{4a^2},
\]
by~\eqref{eq:m-choice}.
\end{proof}

\subsection{Realizing the flat block on the models}\label{sec:realization}

In this section we preserve the exact coefficient values and the
support bound while obtaining
\begin{equation}\label{eq:extension}
 \norm{P+Q}_1\le\norm P_1+\eta,
 \qquad
 \HH(P+Q)\le\theta\HH(P)+\eta
\end{equation}
on each group in~\eqref{eq:models}. Each construction works for every
sufficiently small $a>0$, once the old data and $\eta$ have been fixed.
We also protect a whole finite set of old frequencies, including those
with coefficient zero.

\subsubsection{Products of finite cyclic groups}

Let $G=(\Z/p\Z)^{\N}$ and suppose that $P$ depends on the first $n$
coordinates. Take $Y=(\Z/p\Z)^n$, set
$E=\wh Y=(\Z/p\Z)^n$, and use $m+1$ fresh coordinates for
$z_1,\ldots,z_m,w$. Their joint distribution is exactly the product
distribution in Lemma~\ref{lem:product}, with $\mathbb U$ the $p$th
roots of unity. When $p=2$, all old characters are real and
$\kappa=1$; when $p\ge3$, use $\kappa=2$.

The resulting polynomial $Q$ has coefficients in $\{0,a\}$ and
satisfies~\eqref{eq:extension} with $D=p^n$ in~\eqref{eq:block-count}.
Its support misses the whole old subgroup
$(\Z/p\Z)^n\subset\Gamma$, because the fresh $w$ coordinate has
nonzero exponent. All new frequencies lie in the subgroup on the
first $n+m+1$ coordinates. No averaging error occurs in this case.

\subsubsection{Torsion duals of unbounded exponent}

Suppose that $\Gamma$ is torsion of unbounded exponent and $F\le\Gamma$
is a finite subgroup. The quotient $\Gamma/F$ also has unbounded exponent.
Indeed, if $N(\Gamma/F)=0$ and $MF=0$, then $MN\Gamma=0$.
Let $\pi_F:G\to\wh F$ be the natural quotient map.

\begin{lemma}[Finite averaging modulo the old subgroup]\label{lem:torsion-averaging}
Given finitely many functions $f_\ell\in C(\wh F\times\T^s)$,
$\tau>0$, and an integer $B_0\ge1$, there are
$\gamma\in\Gamma$ and positive integers $N_1,\ldots,N_s$ such that
\begin{equation}\label{eq:torsion-average}
 \left|\int_G f_\ell\bigl(\pi_F(x),
      (N_1\gamma)(x),\ldots,(N_s\gamma)(x)\bigr)\,dm_G(x)
       -\int_{\wh F\times\T^s}f_\ell\right|<\tau
\end{equation}
for every $\ell$, and
\begin{equation}\label{eq:torsion-no-relations}
 h+\sum_{j=1}^s v_jN_j\gamma\ne0
 \quad\text{if }h\in F,\quad
 0\ne v\in\Z^s,\quad |v_j|\le B_0.
\end{equation}
\end{lemma}

\begin{proof}
Approximate each $f_\ell$ uniformly, to accuracy less than $\tau/2$,
by a finite linear combination of characters $h(y)z^v$, with $h\in F$.
Take an integer $B\ge\max(B_0,2)$ bounding every occurring $|v_j|$.
Choose an integer $b>2B+1$ and set $N_j=b^{j-1}$.
If $v\ne0$ and $|v_j|\le B$, the contribution corresponding to the
largest index $j$ with $v_j\ne0$ in
\[
 k(v)=\sum_{j=1}^s v_jb^{j-1}
\]
has modulus greater than the sum of all lower-index contributions.
Hence $k(v)\ne0$, and $|k(v)|\le B\sum_jN_j$.
Choose $\gamma$ whose order modulo $F$ exceeds $B\sum_jN_j$.
Then $k(v)\gamma\notin F$ for all such nonzero $v$.
This proves~\eqref{eq:torsion-no-relations}.

After substitution, a monomial with $v\ne0$ is a nontrivial character
of $G$ and has integral zero. When $v=0$, the monomial has the same
integral as $h$ on $\wh F$, since $\pi_F$ pushes Haar measure to
Haar measure. Thus the integrals of the approximating polynomials
agree exactly. The two uniform approximation errors prove
\eqref{eq:torsion-average}.
\end{proof}

Suppose now that $\supp\wh P\subset F$. Apply
Lemma~\ref{lem:product} with $Y=\wh F$, $E=F$, $D=|F|$,
$\mathbb U=\T$ and $\kappa=2$, allowing error $\eta/2$.
After $d,a,m$ have been chosen, apply
Lemma~\ref{lem:torsion-averaging} to the two continuous functions
\[
 |P+\mathcal Q|,\qquad \sqrt{|P+\mathcal Q|},
\]
with $s=m+1$, $\tau<\eta/2$ and $B_0=2$.
Substitute $z_j=(N_j\gamma)(g)$ and $w=(L\gamma)(g)$, where $g\in G$
and $L=N_{m+1}$. The unreflected copy of the new Fourier support consists
of frequencies of the form
\begin{equation}\label{eq:torsion-frequencies}
 h+\left(L+\sum_{j=1}^m x_jN_j\right)\gamma,
 \qquad h\in F,\quad x\in\{0,1\}^m;
\end{equation}
the reflected copy consists of their negatives.

Condition~\eqref{eq:torsion-no-relations} verifies every required
separation. A frequency in~\eqref{eq:torsion-frequencies} cannot
belong to $F$, since its $L$ coefficient is one. If two frequencies
in the same copy agree and their bit vectors differ, their difference
has a nonzero vector of new exponents in $\{-1,0,1\}^{m+1}$,
which is impossible. If their bit vectors agree, equality forces
the two old characters to agree. Finally, a frequency from one copy
cannot equal a frequency from the other: the resulting relation
has coefficient two at $L$ and coefficients in $\{0,1,2\}$ elsewhere.
This is the reason for requiring $B_0=2$.

Thus the realized coefficients remain exactly $0$ or $a$.
The support misses $F$, lies in the finite subgroup $\langle F,\gamma\rangle$,
and satisfies~\eqref{eq:block-count}. The averaging estimate proves
\eqref{eq:extension}. No direct-sum assumption about
$F$ and $\langle\gamma\rangle$ has been used.

\subsubsection{The circle}

Write $e_k(t)=\exp(2\pi ikt)$ on $\T=\mathbb R/\Z$.
Suppose that
\[
 P=\sum_{|k|\le F_0}p_ke_k,
\]
and that all frequencies in $[-H,H]$ are to be protected, where
$H\ge F_0$ is an integer. Use Lemma~\ref{lem:product} with
$Y=\T$, $E=[-F_0,F_0]\cap\Z$, $D=2F_0+1$,
$\mathbb U=\T$ and $\kappa=2$, allowing error $\eta/2$.

Here is the finite averaging argument with the additional support
conditions included. Uniformly approximate the two continuous
functions $|P+\mathcal Q|$ and $\sqrt{|P+\mathcal Q|}$ on
$\T\times\T^m\times\T$ by trigonometric polynomials.
Let $B\ge1$ bound the absolute values of all their exponents,
including the exponent in the old $\T$ coordinate.
Choose $N_1,\ldots,N_m,L$ successively so that, with
$N_0=1$ and with $L$ regarded as $N_{m+1}$,
\begin{equation}\label{eq:circle-lacunarity}
 N_j>2B\sum_{i=0}^{j-1}N_i\quad(1\le j\le m+1).
\end{equation}
Also require
\begin{equation}\label{eq:circle-separation}
 N_j>2F_0+\sum_{i<j}N_i\quad(1\le j\le m),
 \qquad
 L>H+F_0+\sum_{i=1}^mN_i.
\end{equation}
All these are lower bounds, so they can be imposed simultaneously.

Substitute $z_j=e_{N_j}(t)$ and $w=e_L(t)$.
In a nonconstant monomial of either approximating polynomial, the
contribution from the largest index with nonzero new exponent dominates
all lower-frequency terms by~\eqref{eq:circle-lacunarity}. If all new exponents vanish, a
nonzero old exponent is already nontrivial. Thus every nonconstant
monomial has integral zero after substitution, exactly as on the
product. Uniform approximation gives an arbitrarily small error in
the two integrals, proving~\eqref{eq:extension}.

The positive-frequency part of the realized Fourier support consists of
the frequencies
\[
 L+k+\sum_{j=1}^m x_jN_j,
 \qquad |k|\le F_0,\quad x_j\in\{0,1\}.
\]
They all exceed $H$ by~\eqref{eq:circle-separation}.
If two bit vectors differ, the largest differing $N_j$ dominates
the lower terms and the possible old difference of size $2F_0$;
therefore the frequencies are distinct. When the vectors agree,
distinct old $k$ give distinct frequencies. The negative-frequency part lies below $-H$ and is disjoint from the
positive-frequency part.
The coefficients are therefore exactly $0$ or $a$, the whole interval
$[-H,H]$ is protected, and~\eqref{eq:block-count} is preserved.

In each of the three model groups, once the old polynomial, the protected
finite set and the tolerance are fixed, the block can be realized for
\emph{every} sufficiently small value of $a$. The subsequent choices of
frequencies impose only lower bounds. This is the quantifier used in the
iteration below.

\subsection{Iteration and stabilization of the coefficients}\label{sec:iteration}

Fix the sequence in Theorem~\ref{thm:main}, put $A=a_1$, and choose
positive numbers $\eta_j$ with
\[
 \sum_{j\ge1}\eta_j\le\eps.
\]
Start with $P_0=A$. We construct real even polynomials
\begin{equation}\label{eq:iteration}
 P_j=P_{j-1}+Q_j
\end{equation}
whose Fourier coefficients all belong to $K$, keeping every
previously protected coefficient, including every protected zero,
unchanged.

On $(\Z/p\Z)^{\N}$ let the protected subgroup $F_{j-1}$ consist of
the characters on all coordinates used so far; start with
$F_0=\{0\}$. In the extension step use all of $F_{j-1}$ as the old
set $E$, and protect the subgroup on the first $n+m+1$ coordinates
afterwards. The number of coordinates strictly increases, so these
subgroups exhaust the dual.

For a countable torsion dual of unbounded exponent, fix an enumeration
$\lambda_1,\lambda_2,\ldots$ of $\Gamma$. Start with $F_0=\{0\}$.
At step $j$ apply the torsion construction with old subgroup
$F_{j-1}$ and then put
\[
 F_j=\langle F_{j-1},\gamma_j,\lambda_j\rangle.
\]
This is finite, contains $\supp\wh P_j$, and the subgroups $F_j$
increase to $\Gamma$. Frequencies added only to enlarge the subgroup
are protected with their current coefficients, usually zero.

On the circle put $L_0=1$. At step $j$ choose an integer $F_0^{(j)}$
with $\supp\wh P_{j-1}\subset[-F_0^{(j)},F_0^{(j)}]$ and
$F_0^{(j)}\le L_{j-1}$. Use the circle construction with
$H=2L_{j-1}$, so
\begin{equation}\label{eq:circle-gap}
 \supp\wh Q_j\ \cap[-2L_{j-1},2L_{j-1}]=\varnothing.
\end{equation}
Then choose $L_j\ge\max(L_{j-1},j,1,\deg P_j)$.
Here $\deg P$ denotes the largest absolute value of a frequency in
its support. The integers $L_j$ tend to infinity.

For all three models the old finite set and its size $D_j$ are fixed
\emph{before} selecting a new coefficient value. In the torsion
models $D_j=|F_{j-1}|$, and on the circle
$D_j=2F_0^{(j)}+1$.
Choose a value $c_j$ from the prescribed sequence so far out that the
extension with error $\eta_j$ is available and
\begin{equation}\label{eq:c-choice}
 0<c_j<\min(A,1,2^{-j},c_{j-1}/2),
 \qquad D_jc_j^{1/j}\le2^{-j},
\end{equation}
where $c_0=A$. This is possible because the extension works for
every sufficiently small positive value and $a_n\to0$.
The new polynomial has
\begin{equation}\label{eq:flat-data}
 \wh Q_j=c_j\one_{E_j},\qquad
 E_j=\supp\wh Q_j,\qquad
 |E_j|\le\frac{A^2D_j}{4c_j^2}.
\end{equation}
The finite sets $E_j$ are pairwise disjoint and miss $0$.
No previously protected coefficient changes, so the coefficients
of $P_j$ still belong to $K$ and lie between zero and $A$.
From~\eqref{eq:extension},
\begin{equation}\label{eq:iteration-bounds}
 \norm{P_j}_1\le A+\sum_{i=1}^j\eta_i\le A+\eps,
 \qquad
 \HH(P_j)\le\theta\HH(P_{j-1})+\eta_j.
\end{equation}
In particular,
\begin{equation}\label{eq:H-limit}
 \HH(P_j)\longrightarrow0.
\end{equation}
Indeed, iteration bounds it by
$\theta^j\sqrt A+\sum_{i=1}^j\theta^{j-i}\eta_i$.
For any fixed initial segment, the factors $\theta^{j-i}$ tend to zero;
the remaining terms are bounded by the corresponding tail of the
convergent series $\sum_i\eta_i$.

Every character is eventually protected. Thus the Fourier
coefficients of $P_j$ stabilize, with value $A$ at $0$, value $c_j$
on $E_j$, and zero elsewhere. The uniform variation bound in
\eqref{eq:iteration-bounds} implies that $\int fP_j\,dm_G$ converges
for every $f\in C(G)$. To verify this directly, first use stabilization
for trigonometric polynomials, and then approximate $f$ uniformly
by such polynomials. The resulting limit is a bounded linear
functional of norm at most $A+\eps$ on $C(G)$.
By the Riesz representation theorem there is a measure $\mu$ with
\begin{equation}\label{eq:weakstar}
 P_jm_G\overset{*}{\longrightarrow}\mu,
 \qquad
 \wh\mu=A\one_{\{0\}}+\sum_{j\ge1}c_j\one_{E_j}.
\end{equation}
This proves convergence of the whole sequence; no assertion about
sequential compactness of a weak-star unit ball is needed.

It follows already that $\wh\mu(0)=A$,
$A\le\norm\mu\le A+\eps$, and $\wh\mu(\Gamma)\subset K$.
Reality and evenness pass to the limit. Since $c_j\to0$ and each
$E_j$ is finite, $\wh\mu\in c_0(\Gamma)$.

We record the structural information which will also be used for
complex sequences. The values $c_j=a_{n_j}$ have strictly increasing
indices $n_j\ge2$ by~\eqref{eq:c-choice}, and
\begin{equation}\label{eq:increment-norm}
 \norm{Q_j}_1\le\norm{P_j}_1+\norm{P_{j-1}}_1
 \le2(A+\eps).
\end{equation}
After Haar lifting from a model quotient, all of
\eqref{eq:iteration}, \eqref{eq:weakstar},
and~\eqref{eq:increment-norm} remain true on the original group:
pull back the polynomials and extend their Fourier transforms by zero.

\subsection{Exact recovery and singularity}\label{sec:singularity}

The small values of $\HH(P_j)$ do not by themselves determine the
Lebesgue decomposition of an arbitrary weak-star limit.
Here singularity follows because the polynomials are exact
convolutions of the limit with suitable kernels.
This recovery-and-vanishing strategy has a classical antecedent in
Katznelson's proof of singularity for separated products on the circle
\cite[Lemma~2 and its proof, p.~719]{Katznelson1965}: suitable partial
products recover the density of the absolutely continuous part, and
their limiting behaviour forces that density to vanish.
We implement this strategy using explicit convolution kernels and
convergence in measure, allowing signed measures and covering all
three model groups. The kernel argument is included in full.

\begin{lemma}[Recovery of the absolutely continuous part]\label{lem:ac-recovery}
Let $G$ be a compact Hausdorff abelian group and suppose that
$K_j\in C(G)$ satisfy
\begin{equation}\label{eq:kernel-assumptions}
 \int_G K_j\,dm_G=1,\qquad
 \sup_j\norm{K_j}_1\le C<\infty,
\end{equation}
and, for every open neighbourhood $V$ of $0$,
\begin{equation}\label{eq:kernel-localization}
 \sup_{x\notin V}|K_j(x)|\longrightarrow0.
\end{equation}
If $\nu=fm_G+\nu_s$ is the Lebesgue decomposition of a finite
complex regular Borel measure, then $K_j*\nu\to f$ in Haar measure.
\end{lemma}

\begin{proof}
Let $g\in C(G)$ and $\alpha>0$. Uniform continuity gives an identity
neighbourhood $V$ such that $|g(x-y)-g(x)|<\alpha$ for all
$x\in G$ and $y\in V$. Since $\int K_j=1$,
\[
 \norm{K_j*g-g}_\infty
 \le C\alpha+2\norm g_\infty\sup_{y\notin V}|K_j(y)|.
\]
Thus $K_j*g\to g$ uniformly. Convolution by $K_j$ has $L^1$
operator norm at most $C$, and $C(G)$ is dense in $L^1(G)$.
It follows that $K_j*f\to f$ in $L^1(G)$ for every $f\in L^1(G)$.

Fix $\alpha>0$. By singularity and regularity, choose a compact set
$E$ with $m_G(E)=0$ and
$\norm{\nu_s-\nu_s|_E}<\alpha$.
Choose an open $U\supset E$ with $m_G(U)<\alpha$.
Compactness of $E$ gives an open symmetric identity neighbourhood
$V$ such that $E+V\subset U$.
For $x\notin U$ and $y\in E$ we have $x-y\notin V$, whence
\[
 \sup_{x\notin U}|K_j*(\nu_s|_E)(x)|
 \le\norm{\nu_s}\sup_{z\notin V}|K_j(z)|\longrightarrow0.
\]
The remainder satisfies
\[
 \norm{K_j*(\nu_s-\nu_s|_E)}_1\le C\alpha.
\]
For each $t>0$, the preceding uniform estimate and Markov's
inequality yield
\[
 \limsup_{j\to\infty}m_G\{|K_j*\nu_s|>t\}
 \le\alpha+\frac{2C\alpha}{t}.
\]
Letting $\alpha\downarrow0$ proves $K_j*\nu_s\to0$ in measure.
Combining the two parts of the decomposition proves the lemma.
\end{proof}

\subsubsection{Finite-subgroup kernels}

In either torsion model, the construction provides increasing finite
subgroups $F_j$ exhausting $\Gamma$. Put
\[
 D_{F_j}=\sum_{\gamma\in F_j}\gamma
        =|F_j|\one_{F_j^\perp},
 \qquad
 F_j^\perp=\{x\in G:\gamma(x)=1\text{ for every }\gamma\in F_j\}.
\]
The equality follows by summing the characters of a finite group.
These kernels are continuous and nonnegative, with integral and
$L^1$ norm one. Since characters separate points,
$\bigcap_jF_j^\perp=\{0\}$.
For any open identity neighbourhood $V$, one has
$F_j^\perp\subset V$ for all sufficiently large $j$.
Otherwise the nested nonempty compact sets
$F_j^\perp\cap(G\setminus V)$ would have nonempty intersection,
a contradiction. Thus $D_{F_j}$ is eventually zero outside $V$, and the sequence
$(D_{F_j})$ satisfies the hypotheses of Lemma~\ref{lem:ac-recovery}.

Since $\wh D_{F_j}=\one_{F_j}$, stabilization of the entire
protected subgroup gives the exact identity
\begin{equation}\label{eq:torsion-exact}
 D_{F_j}*\mu=P_j.
\end{equation}
In particular, coefficients which were zero when $F_j$ was protected
contribute zero on both sides.

\subsubsection{Circle kernels}

Let
\[
 \mathsf F_N(t)=\frac1N\left|\sum_{k=0}^{N-1}e_k(t)\right|^2,
 \qquad V_N=2\mathsf F_{2N}-\mathsf F_N.
\]
The kernel $V_N$ is the classical \emph{de la Vall\'ee Poussin kernel};
see \cite[Section~I.2.13]{Katznelson2004}, with the indexing adjusted to
our convention for the Fej\'er kernels.
The Fej\'er kernels have integral and $L^1$ norm one, so
$\int V_N=1$ and $\norm{V_N}_1\le3$.
Their Fourier coefficients give
\begin{equation}\label{eq:V-multiplier}
 \wh V_N(n)=1\quad(|n|\le N),\qquad
 \wh V_N(n)=0\quad(|n|\ge2N).
\end{equation}
For $0<\delta\le1/2$, the geometric-sum formula implies
\[
 \sup_{\norm{t}_{\T}\ge\delta}|V_N(t)|
 \le\frac{2}{N\sin^2(\pi\delta)}\longrightarrow0.
\]
Here $\norm{t}_{\T}$ is the distance from $t$ to zero on the circle.
Thus the sequence $(V_{L_j})$ satisfies the hypotheses of
Lemma~\ref{lem:ac-recovery}.

The Fourier support of $P_j$ lies in $[-L_j,L_j]$.
Every later block misses $[-2L_j,2L_j]$ by
\eqref{eq:circle-gap} and the monotonicity of $L_j$.
Together with~\eqref{eq:V-multiplier}, this proves
\begin{equation}\label{eq:circle-exact}
 V_{L_j}*\mu=P_j.
\end{equation}
Indeed, the transform of the left-hand side agrees with that of
$P_j$ on $[-L_j,L_j]$, is zero in the intervening protected band,
and vanishes beyond $2L_j$.

\subsubsection{Conclusion of the singularity argument}

Write $\mu=fm_G+\mu_s$ on any model group.
Equations~\eqref{eq:torsion-exact} and~\eqref{eq:circle-exact},
and Lemma~\ref{lem:ac-recovery}, show that $P_j\to f$ in measure.
On the other hand, for every $t>0$,
\[
 m_G\{|P_j|>t\}\le t^{-1/2}\HH(P_j)\longrightarrow0
\]
by~\eqref{eq:H-limit}.
Uniqueness of the limit in measure gives $f=0$ almost everywhere.
Hence $\mu\perp m_G$ on each model, and Haar lifting preserves
this conclusion on the original group.

\subsection{Summability, strong continuity, and spectra}\label{sec:completion}

\subsubsection{Summability of Fourier coefficients and convolution powers}

Fix $r>2$. For all sufficiently large $j$ we have $1/j\le r-2$.
Since $c_j<1$, the block bound and the choice of $c_j$ give
\begin{equation}\label{eq:ell-r}
 \sum_{\gamma\in E_j}|\wh\mu(\gamma)|^r
 =|E_j|c_j^r
 \le\frac{A^2}{4}D_jc_j^{r-2}
 \le\frac{A^2}{4}D_jc_j^{1/j}
 \le\frac{A^2}{4}2^{-j}.
\end{equation}
There are only finitely many earlier blocks, each finite.
This proves $\wh\mu\in\ell^r(\Gamma)$ for every $r>2$.
In particular $(\wh\mu)^2\in\ell^2(\Gamma)$.
Plancherel's theorem gives $g\in L^2(G)$ with
$\wh g=(\wh\mu)^2$, which is also the Fourier transform of
$\mu*\mu$. Fourier uniqueness yields $\mu*\mu=gm_G$.

More generally, let $2\le p<\infty$ and put
$s=p/(p-1)\in(1,2]$. Since $2s>2$, we have
$\wh\mu\in\ell^{2s}(\Gamma)$ and hence
$(\wh\mu)^2\in\ell^s(\Gamma)$. Apply the Hausdorff--Young inequality
to the inverse Fourier transform on the discrete group $\Gamma$, whose
dual is $G$. It gives $g_p\in L^p(G)$ with
$\wh{g_p}=(\wh\mu)^2$; Fourier uniqueness
then identifies $\mu*\mu$ with $g_pm_G$. For $1\le p<2$, the conclusion
follows from $L^2(G)\subset L^p(G)$, since Haar measure is normalized.
Thus
\begin{equation}\label{eq:all-Lp}
 \mu*\mu\in\bigcap_{1\le p<\infty}L^p(G).
\end{equation}
Finally, $\wh\mu\in\ell^3(\Gamma)$, so
$(\wh\mu)^3\in\ell^1(\Gamma)$. Its inverse Fourier series converges
absolutely and uniformly to a continuous function $h$, and Fourier
uniqueness gives $\mu^{*3}=hm_G$. All these conclusions are preserved by
Lemma~\ref{lem:lift}.

\subsubsection{Rajchman measures are strongly continuous}

We supply the needed band property, including the total variation step,
for arbitrary compact abelian groups, including nonmetrizable ones.

\begin{lemma}\label{lem:Rajchman-band}
If $\nu\in M_0(G)$ and $f\in L^1(|\nu|)$, then $f\nu\in M_0(G)$.
In particular, $|\nu|\in M_0(G)$.
Every Rajchman measure is strongly continuous.
\end{lemma}

\begin{proof}
For a trigonometric polynomial $p$, the transform of $p\nu$ is a
finite linear combination of translates of $\wh\nu$, hence belongs
to $c_0(\Gamma)$. Trigonometric polynomials are dense in
$L^1(|\nu|)$: first approximate by continuous functions using
regularity, then use Stone--Weierstrass.
Moreover,
\[
 \sup_{\gamma\in\Gamma}
 |\wh{f\nu}(\gamma)-\wh{p\nu}(\gamma)|
 \le\int_G|f-p|\,d|\nu|.
\]
Since $c_0(\Gamma)$ is uniformly closed, $f\nu$ is Rajchman.
Writing $\nu=u|\nu|$ with $|u|=1$ almost everywhere and taking
$f=\overline u$ proves the variation assertion.

A positive Rajchman measure $\rho$ on a compact abelian
group cannot have an atom of mass $c>0$ at a point $x$.
Otherwise $c^{-1}\one_{\{x\}}\rho=\delta_x$ would be Rajchman,
although its Fourier coefficients have modulus one on the infinite
dual group.

Now let $H\subset G$ be a closed subgroup of infinite index, and
let $\pi:G\to G/H$ be the quotient map.
The positive measure $\pi_*|\nu|$ is Rajchman, since its Fourier
transform is the restriction of $\wh{|\nu|}$ to the subgroup
$H^\perp\subset\Gamma$.
The group $G/H$ is infinite, so this pushforward has no atoms.
Consequently
\[
 |\nu|(x+H)=(\pi_*|\nu|)(\{\pi(x)\})=0,
\]
as required.
\end{proof}

\subsubsection{Natural spectrum}

The following elementary implication explains the compatibility with
the Wiener--Pitt results mentioned in the introduction. It belongs to
the natural-spectrum framework for Rajchman measures developed by
Zafran \cite[proof of Proposition~2.9, p.~615]{Zafran1973}; the short proof
is included for completeness.

\begin{lemma}\label{lem:natural-spectrum}
If $\nu*\nu\in L^1(G)$, then $\nu$ has a natural
spectrum.
\end{lemma}

\begin{proof}
The Riemann--Lebesgue lemma gives
$(\wh\nu)^2=\wh{\nu*\nu}\in c_0(\Gamma)$, and therefore
$\wh\nu\in c_0(\Gamma)$. Thus $\nu\in M_0(G)$.
Let $\varphi$ be a nonzero multiplicative linear functional on $M(G)$.
Its restriction to the ideal $L^1(G)$ is either zero or a Fourier
evaluation at a character $\gamma\in\Gamma$; this is the usual
description of the maximal ideal space of the group algebra
\cite{GrahamMcGehee1979,Rudin1962}.
In the first case
\[
 \varphi(\nu)^2=\varphi(\nu*\nu)=0.
\]
In the second case choose $h\in L^1(G)$ with $\wh h(\gamma)\ne0$.
Since $\nu*h\in L^1(G)$,
\[
 \varphi(\nu)\wh h(\gamma)
 =\varphi(\nu*h)
 =\wh\nu(\gamma)\wh h(\gamma),
\]
so $\varphi(\nu)=\wh\nu(\gamma)$.
The spectrum is therefore contained in $\wh\nu(\Gamma)\cup\{0\}$.
The Rajchman condition and the infinitude of $\Gamma$ imply
$0\in\overline{\wh\nu(\Gamma)}$.
Conversely, every Fourier evaluation is a multiplicative functional
on $M(G)$, and the spectrum is closed. This gives equality with the closure
of the Fourier--Stieltjes range.
\end{proof}

\begin{proof}[Completion of the proof of Theorem~\ref{thm:main}]
Sections~\ref{sec:selectors}--\ref{sec:iteration} give a real even
Rajchman measure on every model group, with the stated range and
norm bounds. Section~\ref{sec:singularity} proves singularity, and
\eqref{eq:ell-r} gives every asserted summability condition and the
claimed regularity of the second and third convolution powers. The reduction in
Section~\ref{sec:reduction} lifts these properties to any
compact abelian group.
Lemma~\ref{lem:Rajchman-band} proves strong continuity, and
Lemma~\ref{lem:natural-spectrum} shows that the spectrum is natural.
All Fourier values lie in $[0,A]$, and the value $A$ occurs at the
trivial character. Thus the spectral radius is exactly $A=a_1$.
\end{proof}

\section{Exact palettes and group structure}\label{sec:exact-palettes}

\subsection{Universal exact realization}

We retain the notation~\eqref{eq:palette}. The value $0$ is part of the
palette and must be attained. This requirement is harmless in the
constructions below because one protected frequency is assigned
coefficient zero from the beginning.

Recall from Lemma~\ref{lem:Rajchman-band} that every Rajchman measure on
a compact abelian group is strongly continuous.

\subsection{Torsion duals: forcing finite tranches at bounded cost}

Let $E\le\Gamma$ be finite and put
\[
 D_E(x)=\sum_{\gamma\in E}\gamma(x).
\]

\begin{lemma}\label{lem:dirichlet-subgroup}
For every finite subgroup $E\le\Gamma$,
\[
 D_E=|E|\one_{E^\perp}\ge0,\qquad
 \norm{D_E}_1=1,\qquad \wh{D_E}=\one_E.
\]
The polynomial $D_E$ is real and even.
\end{lemma}

\begin{proof}
If $x\in E^\perp$, every term equals one. Otherwise choose
$\eta\in E$ with $\eta(x)\ne1$; translation of the index set gives
$D_E(x)=\eta(x)D_E(x)$, hence $D_E(x)=0$. Since
$m_G(E^\perp)=|E|^{-1}$, the norm identity follows. The Fourier identity
is the character orthogonality relation.
\end{proof}

\begin{lemma}[Finite-tranche forcing]\label{lem:finite-tranche}
Assume that $\Gamma$ is an infinite torsion group. Let $F\le\Gamma$ be
finite and $p\le q$. There is a real even trigonometric polynomial $S$
such that
\begin{gather*}
 \wh S|_F=0,\qquad
 \wh S(\Gamma)\subset\{0,a_p,\ldots,a_q\},\\
 \{a_p,\ldots,a_q\}\subset\wh S(\Gamma),
 \qquad \norm S_1\le2a_p.
\end{gather*}
Every displayed nonzero value occurs only finitely many times.
\end{lemma}

\begin{proof}
Build a chain of proper finite subgroup extensions
\[
 F=E_0\subsetneq E_1\subsetneq\cdots\subsetneq E_m,
 \qquad m=q-p+1,
\]
and set
\[
 S=\sum_{r=1}^{m}a_{p+r-1}(D_{E_r}-D_{E_{r-1}}).
\]
Its transform equals $a_{p+r-1}$ on $E_r\setminus E_{r-1}$ and is zero
on $F$ and outside $E_m$. Moreover,
\[
 S=R-a_pD_F,
\]
where
\[
 R=a_qD_{E_m}+\sum_{r=1}^{m-1}
 (a_{p+r-1}-a_{p+r})D_{E_r}\ge0,
 \qquad \int_G R\,dm_G=a_p.
\]
Thus $\norm S_1\le2a_p$.
\end{proof}

The point is that the cost is governed by the largest value in the
tranche and not by the number of prescribed values. In particular,
\begin{equation}\label{eq:H-forcing}
 \HH(P+S)\le\HH(P)+\sqrt{\norm S_1}
 \le\HH(P)+\sqrt{2a_p}.
\end{equation}

The torsion case of Lemma~\ref{lem:algebra} will be used repeatedly:
every infinite torsion abelian group contains either a countable torsion
subgroup of unbounded exponent or a subgroup isomorphic to
$\bigoplus_{n\ge1}\Z/p\Z$ for some prime $p$.

\begin{theorem}[Exact palettes on torsion duals]
\label{thm:exact-torsion}
Let $G$ be a compact abelian group with torsion dual. For every
positive null palette $K$ and every $\eps>0$ there is a real even measure
$\mu$ such that
\begin{gather}
 \mu\perp m_G,
 \qquad \mu\in M_0(G),
 \qquad \wh\mu(\Gamma)=K,\label{eq:exact-torsion-range}\\
 \mu(G)=a_1,
 \qquad a_1\le\norm\mu\le a_1+\eps.
 \label{eq:exact-torsion-norm}
\end{gather}
Consequently, $\mu$ is atomless and strongly continuous.
\end{theorem}

\begin{proof}
By Lemmas~\ref{lem:algebra} and~\ref{lem:lift}, it suffices
to work on one of the two countable model duals used in
Section~\ref{sec:realization}. Enumerate the model dual as
$\Gamma=\{\lambda_1,\lambda_2,\ldots\}$.

Choose indices $2\le p_1<p_2<\cdots$ so rapidly that
\[
 a_{p_j}\le\min\{\eps2^{-j-2},16^{-j}\},
\]
and choose $\eta_j>0$ with $\sum_j\eta_j\le\eps/2$. The first tranche
is $a_1,\ldots,a_{p_1-1}$; the $j$th subsequent tranche is
$a_{p_j},\ldots,a_{p_{j+1}-1}$.

Let $m=p_1-1$ and choose a strictly increasing chain of finite subgroups
\[
 \{0\}=E_0\subsetneq E_1\subsetneq\cdots\subsetneq E_{m-1},
\]
with just $E_0$ when $m=1$. Set
\[
 P_0=a_mD_{E_{m-1}}+
     \sum_{r=1}^{m-1}(a_r-a_{r+1})D_{E_{r-1}}.
\]
Lemma~\ref{lem:dirichlet-subgroup} shows that $P_0$ is nonnegative,
real and even, and
\[
 \norm{P_0}_1=\int_G P_0\,dm_G
 =a_m+\sum_{r=1}^{m-1}(a_r-a_{r+1})=a_1.
\]
Its transform equals $a_1$ on $E_0$, equals $a_{r+1}$ on
$E_r\setminus E_{r-1}$ for $1\le r\le m-1$, and vanishes outside
$E_{m-1}$. Thus every value in the first tranche is attained. Choose
$\xi\in\Gamma\setminus E_{m-1}$ and a finite subgroup $F_0$ containing
$E_{m-1}\cup\{\xi\}$. The coefficient $\wh P_0(\xi)=0$ is protected
from the beginning. Set also $c_0=a_1$.

Suppose that $P_{j-1}$ and a finite protected subgroup $F_{j-1}$ have
been constructed. Apply Lemma~\ref{lem:finite-tranche} to obtain $S_j$
which realizes the next tranche outside $F_{j-1}$ and satisfies
$\norm{S_j}_1\le2a_{p_j}$. Put $R_j=P_{j-1}+S_j$ and enlarge the
protected subgroup so that it contains $\supp\wh R_j$ and $\lambda_j$.
Apply the flat-block realization from Section~\ref{sec:realization} to
$R_j$, with error $\eta_j$, choosing the new coefficient $c_j$ from a
sufficiently far tail of the palette and requiring
\[
 0<c_j<\min\{a_{p_j},2^{-j},c_{j-1}/2\}.
\]
This gives a block $Q_j$ with coefficients $0$ or $c_j$, supported
outside the protected subgroup, and
\begin{align}
 \norm{P_j}_1&\le\norm{P_{j-1}}_1+2a_{p_j}+\eta_j,
 \label{eq:exact-norm-rec}\\
 \HH(P_j)&\le\theta\HH(P_{j-1})
  +\theta\sqrt{2a_{p_j}}+\eta_j,
 \qquad P_j=R_j+Q_j,\label{eq:exact-H-rec}
\end{align}
where $0<\theta<1$ is absolute. Enlarge once more to a finite subgroup
$F_j$ containing the new support. The groups can be chosen increasing
and exhausting $\Gamma$.

The norm recursion gives
\[
 \sup_j\norm{P_j}_1\le a_1+2\sum_j a_{p_j}+\sum_j\eta_j
 \le a_1+\eps.
\]
The inhomogeneous term in~\eqref{eq:exact-H-rec} tends to zero, hence
$\HH(P_j)\to0$. Coefficients on $F_j$ are protected forever. Therefore
a weak-star limit $\mu$ has exactly the stabilized coefficients. All
palette values occur because the tranches cover every index; no other
nonzero value is introduced, and the protected character $\xi$ retains
coefficient zero. This proves $\wh\mu(\Gamma)=K$ and $\mu(G)=a_1$.

The size of values introduced after stage $j$ tends to zero, while each
stage introduces only finitely many frequencies, so $\mu\in M_0(G)$.
For exact recovery, the subgroup kernels satisfy
\[
 \mu*D_{F_j}=P_j.
\]
They form a contractive approximate identity in $L^1(G)$ and localize at
the identity. Lemma~\ref{lem:ac-recovery} and $\HH(P_j)\to0$ show that
the absolutely continuous part of $\mu$ vanishes. Finally, Haar lifting
back to the original group preserves the norm, mass, exact range,
Rajchman property, and pure singularity.
\end{proof}

\begin{corollary}\label{cor:infinite-torsion}
If $\Tor\Gamma$ is infinite, the conclusion of
Theorem~\ref{thm:exact-torsion} holds on $G$.
\end{corollary}

\begin{proof}
The inclusion $\Tor\Gamma\le\Gamma$ corresponds to a quotient of $G$
whose dual is $\Tor\Gamma$. Apply Theorem~\ref{thm:exact-torsion} on the
quotient and lift by Lemma~\ref{lem:lift}.
\end{proof}

\subsection{An explicit telescoping construction from an infinite-order character}

\begin{theorem}\label{thm:telescope}
Assume that $\Gamma$ contains an element $\gamma$ of infinite order.
Let
\[
 H_j=\ker(2^j\gamma)\quad(j\ge0),\qquad
 d_j=a_j-a_{j+1}\quad(j\ge1).
\]
Then
\begin{equation}\label{eq:telescope-measure}
 \mu=a_1m_{H_0}-\sum_{j\ge1}d_jm_{H_j}
\end{equation}
is a real even purely singular measure of norm at most $2a_1$, mass
zero, and exact Fourier--Stieltjes range $K$. If
$[\Gamma:\langle\gamma\rangle]=\infty$, the measure is atomless.
\end{theorem}

\begin{proof}
The series converges absolutely in variation since $\sum_jd_j=a_1$.
Every $H_j$ has Haar measure zero in $G$, so the total variation is
concentrated on a countable union of Haar-null sets. Moreover,
\[
 \wh{m_{H_j}}=\one_{2^j\langle\gamma\rangle}.
\]
Thus the transform vanishes off $\langle\gamma\rangle$. For $k\ne0$, let
\[
 v_2(k)=\max\{j\ge0:2^j\text{ divides }|k|\};
\]
the sum below is understood to be empty when $v_2(k)=0$. Then
\[
 \wh\mu(k\gamma)=a_1-\sum_{j=1}^{v_2(k)}d_j
 =a_{v_2(k)+1},
\]
while $\wh\mu(0)=0$. Every value occurs. If the quotient
$\Gamma/\langle\gamma\rangle$ is infinite, then $H_0$ and all $H_j$
are infinite compact groups, so their Haar measures are atomless.
\end{proof}

For $G=\T$, formula~\eqref{eq:telescope-measure} becomes
\[
 a_1\delta_0-\sum_{j\ge1}(a_j-a_{j+1})m_{C_{2^j}},
\]
where $C_{2^j}$ is the subgroup of $2^j$th roots of unity. Hence every
palette has an atomic purely singular exact realization on the circle,
even though atomless universal realization fails.

\begin{proposition}\label{prop:telescope-mass}
Every palette also has an atomic real even exact realization $\eta$ on
$\T$ satisfying $\eta(\T)=a_1$ and $\norm\eta\le3a_1$.
If $\gamma$ is an infinite-order character on a general $G$ and
$[\Gamma:\langle\gamma\rangle]=\infty$, the Haar lift of $\eta$ along the
quotient map induced by $\gamma$ is atomless.
\end{proposition}

\begin{proof}
With $e_n(t)=\exp(2\pi int)$ on $\T=\mathbb R/\Z$, put
\[
 \eta=a_1\delta_0-\sum_{j\ge1}d_j(e_1+e_{-1})m_{C_{2^{j+1}}}.
\]
The series converges in variation, is real and even, and has mass $a_1$.
Its transform equals $a_1$ off the nested sets of residue classes
$\{n:n\equiv\pm1\pmod{2^{j+1}}\}$, equals $a_{j+1}$ on the difference
between two successive nested sets of residue classes, and equals zero at
$\pm1$. This gives exact
range and the norm estimate. The lifting assertion follows from
Lemma~\ref{lem:lift}.
\end{proof}

\subsection{Gap obstructions}

The negative part starts from the classical quantitative gap theorem
\cite[Theorem~1.4.1, pp.~18--19]{GrahamMcGehee1979}.
The circle case goes back to de Leeuw and Katznelson
\cite[Lemma~2]{DeLeeuwKatznelson1970}. Ramsey \cite{Ramsey1977}
treated compact abelian groups whose dual has finite torsion subgroup,
and Ramsey and Wells \cite{RamseyWells1977} treated arbitrary compact
abelian groups. For the attribution of the two proofs presented in the
book, see \cite[p.~27]{GrahamMcGehee1979}.

For the cardinality estimate one can do substantially better than the
bound obtained by solving the numerical condition in that theorem. Once
the large-coefficient set is known to be finite, the full Littlewood
theorem of McGehee--Pigno--Smith and Konyagin
\cite{McGeheePignoSmith1981,Konyagin1982}, combined with the invariant
uniform approximation property of Bo\.{z}ejko--Pe{\l}czy\'nski and Bourgain
\cite{BozejkoPelczynski1979,Bourgain1987}, gives an estimate exponential
in the norm of the measure. The circle version of this argument is
\cite[Lemma~8]{OW2014}; we include the short extension to torsion-free
dual groups.

\begin{theorem}[Quantitative gap principle]\label{thm:gap-principle}
There are an absolute constant $C>0$ and a nonincreasing function
\[
 \eta:[1,\infty)\longrightarrow(0,1)
\]
with the following property. Let $G$ have torsion-free dual, let
$K\ge1$, and let $\sigma\in M(G)$ be strongly continuous with
$\norm\sigma\le K$. Put
\[
 Q=\{\gamma\in\Gamma:|\wh\sigma(\gamma)|\ge1\}.
\]
If
\[
 |\wh\sigma(\gamma)|\le\eta(K)\qquad(\gamma\notin Q),
\]
then
\[
 |Q|\le\exp(CK).
\]
For $G=\T$, strong continuity is equivalent to atomlessness.
\end{theorem}

\begin{proof}
Let $L>0$ be an absolute Littlewood constant: every trigonometric
polynomial with $N$ nonzero coefficients of modulus at least one has
$L^1$ norm at least $L\log N$. The result of
McGehee--Pigno--Smith and Konyagin gives this on $\T$. The same estimate
holds whenever the dual is torsion-free. Indeed, the finite support
generates a free abelian group $H\cong\Z^s$, so the polynomial factors
through $G/H^\perp\cong\T^s$. Restricting it to
\[
 t\longmapsto(t,Mt,\ldots,M^{s-1}t)
\]
and letting $M\to\infty$ preserves its $L^1$ norm in the limit, while for
large $M$ the resulting one-dimensional frequencies are distinct.

Let $\alpha>1$ be an absolute constant in the invariant uniform
approximation theorem: for every finite $E\subset\Gamma$ there is a
trigonometric polynomial $\Theta$ such that
\[
 \wh\Theta=1\ \hbox{ on }E,\qquad
 \norm\Theta_1<2,\qquad
 |\supp\wh\Theta|\le\alpha^{2|E|}.
\]
For $K\ge1$ set
\[
 d(K)=\left\lceil\exp(4K/L)\right\rceil,
 \qquad r(K)=\left\lceil16K^2\right\rceil+3,
\]
and
\begin{equation}\label{eq:eta-gap}
 \eta(K)=\min\left\{e^{-r(K)},\frac{1}{4\alpha^{2d(K)}}\right\}.
\end{equation}
Both terms on the right-hand side are nonincreasing in $K$, so $\eta$
has the asserted monotonicity.

The first term in~\eqref{eq:eta-gap}, together with
$K<\sqrt{r(K)}/4$, allows us to apply
\cite[Theorem~1.4.1\textup{(A)}]{GrahamMcGehee1979}; hence $Q$ is finite.
Let $f$ be the trigonometric polynomial defined by
\[
 \wh f=\wh\sigma\one_Q
\]
and put $\nu=\sigma-fm_G$. Then every nonzero coefficient of $f$ has
modulus at least one, while
\[
 \sup_{\gamma\in\Gamma}|\wh\nu(\gamma)|\le\eta(K),
 \qquad \norm{fm_G+\nu}=\norm\sigma\le K.
\]

Suppose that $|Q|\ge d(K)$. A dyadic-coset selection gives a coset
$A\subset\Gamma$ of a subgroup such that
\[
 d(K)\le|Q\cap A|<2d(K).
\]
For completeness, if $H$ is generated by $Q$, choose a homomorphism
$q:H\to\Z$ which is injective on $Q$ and successively select a residue
class of $q(Q)$ modulo $2^j$ containing at least $d(K)$ points, stopping
as soon as it contains fewer than $2d(K)$. The indicator of the
corresponding coset of $q^{-1}(2^j\Z)$ is the Fourier--Stieltjes transform
of a measure of norm one. Convolving with this measure therefore produces
$f_1$ and $\nu_1$ satisfying
\[
 d(K)\le|\supp\wh f_1|<2d(K),
 \qquad \norm{f_1m_G+\nu_1}\le K,
 \qquad \sup_\gamma|\wh{\nu_1}(\gamma)|\le\eta(K).
\]
The Littlewood estimate gives $\norm{f_1}_1\ge4K$. Apply the invariant
uniform approximation theorem to $E=\supp\wh f_1$. Since
$|E|<2d(K)$, the resulting polynomial $\Theta$ satisfies
$|\supp\wh\Theta|\le\alpha^{4d(K)}$, and hence Parseval's identity gives
\[
 \norm{\nu_1*\Theta}_1
 \le\norm{\nu_1*\Theta}_2
 \le2\eta(K)\alpha^{2d(K)}\le\frac12.
\]
On the other hand,
\[
 2K>\norm{(f_1m_G+\nu_1)*\Theta}
 =\norm{f_1+\nu_1*\Theta}_1
 \ge4K-\frac12,
\]
a contradiction. Thus $|Q|<d(K)$. Enlarging an absolute constant gives
$d(K)\le\exp(CK)$ for all $K\ge1$.
\end{proof}

\begin{proposition}\label{prop:circle-counterpalette}
There exists a positive null palette $K_*$ such that every atomless
measure $\mu$ on $\T$ with $\wh\mu(\Z)\subset K_*$ omits infinitely
many nonzero values of $K_*$. In particular, no atomless measure exactly
realizes $K_*$.
\end{proposition}

\begin{proof}
Let $C$ and $\eta$ be supplied by
Theorem~\ref{thm:gap-principle}. Choose recursively $b_1=1/2$ and
\[
 K_j=\frac{j}{b_j},\qquad
 N_j=\left\lceil\exp(CK_j)\right\rceil,
 \qquad b_{j+1}=\tfrac14b_j\eta(K_j).
\]
Let $C_j$ contain $N_j+1$ distinct points of $[b_j,2b_j]$ and put
$K_*=\{0\}\cup\bigcup_jC_j$. Enumerate the nonzero elements of $K_*$ in
decreasing order; this makes $K_*$ a positive null palette. If
$j\ge\norm\mu$, apply Theorem~\ref{thm:gap-principle} to
$\sigma_j=\mu/b_j$. The coefficients outside
$Q_j=\{n:|\wh\mu(n)|\ge b_j\}$ have modulus at most
$2b_{j+1}/b_j\le\eta(K_j)$ after normalization, and hence
$|Q_j|\le N_j$. Since $C_j$ has $N_j+1$ distinct values, at least one
value in each sufficiently late cluster is omitted.
\end{proof}

We need a version uniform over compact groups with finitely many
components.

\begin{lemma}\label{lem:components-split}
Let $G_0$ be the connected component of the identity. Then
\begin{equation}\label{eq:component-duality}
 \wh{G/G_0}=\Tor\Gamma.
\end{equation}
If $\Tor\Gamma$ is finite, then $F=G/G_0$ is finite and
$G\cong G_0\times F$ topologically.
\end{lemma}

\begin{proof}
A character is trivial on $G_0$ exactly when its image is finite,
which proves~\eqref{eq:component-duality}. Thus $G/G_0$ has torsion dual
and is profinite. If the torsion subgroup is finite, Pontryagin duality
makes $F=G/G_0$ finite. The compact connected
abelian group $G_0$ is divisible: multiplication by every positive
integer is onto, since the dual multiplication map on the torsion-free
group $\wh{G_0}$ is injective. Divisible abelian groups are injective
$\Z$-modules, so the finite extension $0\to G_0\to G\to F\to0$
splits algebraically. A section from the
finite discrete group $F$ into $G$ is automatically continuous. Its
image meets $G_0$ only at zero, and the resulting map
$G_0\times F\to G$ is a continuous bijection from a compact group to a
Hausdorff group, hence a topological isomorphism.
\end{proof}

\begin{lemma}\label{lem:partial-fourier}
Let $G=L\times F$ with $F$ finite, and let $\mu\in M(G)$. For $f\in F$,
write $\mu_f$ for the restriction
of $\mu$ to $L\times\{f\}$, identified with a measure on $L$, and for
$\chi\in\wh F$ put
\[
 \nu_\chi=\sum_{f\in F}\overline{\chi(f)}\mu_f.
\]
Then
\[
 \wh{\nu_\chi}(\lambda)=\wh\mu(\lambda,\chi),
 \qquad \norm{\nu_\chi}\le\norm\mu.
\]
Atomlessness and strong continuity of $\mu$ pass to every $\nu_\chi$.
\end{lemma}

\begin{proof}
The Fourier identity is immediate from the definition. The components
are mutually singular, so $\sum_f\norm{\mu_f}=\norm\mu$ and the norm
bound follows. The variation of $\nu_\chi$ is dominated by
$\sum_f|\mu_f|$, which also gives the continuity assertions.
\end{proof}

\begin{theorem}[A universal counterpalette]\label{thm:universal-counterpalette}
There is one positive null palette $K_*$ such that, for every compact
abelian group $G$ with finite $\Tor\wh G$, every strongly continuous
measure whose Fourier--Stieltjes range is contained in $K_*$ omits
infinitely many values
of $K_*$. Hence no such measure exactly realizes $K_*$.
\end{theorem}

\begin{proof}
Let $C$ and $\eta$ be supplied by
Theorem~\ref{thm:gap-principle}. Choose recursively $b_1=1/2$ and set
\[
 K_j=\frac{j}{b_j},\qquad
 B_j=\left\lceil\exp(CK_j)\right\rceil,
 \qquad b_{j+1}=\tfrac14b_j\eta(K_j).
\]
Choose $C_j\subset[b_j,2b_j]$ with $jB_j+1$ elements and put
$K_*=\{0\}\cup\bigcup_jC_j$.

Write $G\cong G_0\times F$. The standing convention and finiteness of $F$
imply that $G_0$ is infinite. For an integer
$j\ge\max\{|F|,\norm\mu\}$, apply Lemma~\ref{lem:partial-fourier} and
Theorem~\ref{thm:gap-principle} to $\nu_\chi/b_j$ on the connected group
$G_0$. For each $\chi\in\wh F$, at most $B_j$ frequencies can carry
coefficients of modulus at least $b_j$. Therefore at most
$|F|B_j\le jB_j$ different values from $C_j$ can occur, whereas
$|C_j|=jB_j+1$. At least one value is omitted in every sufficiently
late cluster.
\end{proof}

\subsection{Classification}

\begin{lemma}\label{lem:virtually-cyclic}
If an abelian group $\Gamma$ contains an element of infinite order,
then exactly one of the following holds:
\begin{enumerate}[label=\textup{(\roman*)}]
\item $\Gamma\cong\Z\oplus A$ for a finite group $A$;
\item every infinite cyclic subgroup has infinite index in $\Gamma$.
\end{enumerate}
\end{lemma}

\begin{proof}
If $\langle\gamma\rangle$ has finite index, then $\Gamma$ is finitely
generated of rank one, hence isomorphic to $\Z\oplus A$ with $A$
finite. The converse follows directly from the structure theorem.
\end{proof}

\begin{lemma}\label{lem:T-finite-continuity}
Every atomless measure on $\T\times F$, with $F$ finite, is strongly
continuous.
\end{lemma}

\begin{proof}
A closed subgroup of infinite index cannot project onto all of $\T$;
its projection is therefore finite, and the subgroup itself is finite.
An atomless measure vanishes on all of its cosets.
\end{proof}

\begin{theorem}\label{thm:classification}
The classification in Theorem~\ref{thm:classification-overview} holds.
In all affirmative cases, the realizing measures may be chosen real and
even.
\end{theorem}

\begin{proof}
For $\mathcal U_{\mathrm s}$, if $\Gamma$ is torsion, use
Theorem~\ref{thm:exact-torsion}; if $\Gamma$ contains an infinite-order
element, use Theorem~\ref{thm:telescope}. This proves part~(a).

For atomless realization, groups $\T\times F$ fail by
Lemma~\ref{lem:T-finite-continuity} and
Theorem~\ref{thm:universal-counterpalette}. Conversely, if $\Gamma$ is
torsion, Theorem~\ref{thm:exact-torsion} applies. Otherwise, unless
$G\cong\T\times F$, Lemma~\ref{lem:virtually-cyclic} supplies an
infinite-order $\gamma$ with infinite-index cyclic subgroup, and
the measure constructed in Theorem~\ref{thm:telescope} is atomless. This
proves part~(b).

If $\Tor\Gamma$ is infinite, Corollary~\ref{cor:infinite-torsion} gives
a Rajchman exact realization with the near-optimal norm estimate; by
Lemma~\ref{lem:Rajchman-band} it is strongly continuous. If
$\Tor\Gamma$ is finite, Theorem~\ref{thm:universal-counterpalette}
rules out universal strongly continuous realization. Identity
\eqref{eq:component-duality} gives the equivalence with
$|G/G_0|=\infty$. This proves part~(c).
\end{proof}

In the following table, $p$ is a prime, $\Z_p$ denotes the additive
group of $p$-adic integers, and $\Z(p^\infty)$ is the Pr\"ufer $p$-group.

\begin{center}
\begin{tabular}{>{\raggedright\arraybackslash}p{0.28\textwidth}
                >{\raggedright\arraybackslash}p{0.28\textwidth}ccc}
\toprule
$G$ & $\wh G$ & $\mathcal U_{\mathrm s}$ &
$\mathcal U_{\mathrm c}$ & $\mathcal U_{\mathrm{sc}},\mathcal U_0$\\
\midrule
$\T$ & $\Z$ & yes & no & no\\
$\T\times F$, $F$ finite & $\Z\oplus\wh F$ & yes & no & no\\
$\T^2$ & $\Z^2$ & yes & yes & no\\
$(\Z/2\Z)^{\N}$ & $\bigoplus_{n\ge1}\Z/2\Z$ & yes & yes & yes\\
$\Z_p$ & $\Z(p^\infty)$ & yes & yes & yes\\
$\T\times(\Z/2\Z)^{\N}$ &
$\Z\oplus\bigoplus_{n\ge1}\Z/2\Z$ & yes & yes & yes\\
$\widehat{\Z[1/p]}$ & $\Z[1/p]$ & yes & yes & no\\
\bottomrule
\end{tabular}
\end{center}

\subsection{Palette-dependent criteria on the circle}

Universal atomless realization fails on $\T$, but the following cost
criterion captures several broad positive classes. For finite
$F\subset K\setminus\{0\}$ define
\[
 \lambda_K(F)=\inf\left\{\norm P_1:
 P\text{ is a trigonometric polynomial},\ \wh P(\Z)\subset K,
 \ F\subset\wh P(\Z)\right\},
\]
put $F_{p,q}=\{a_p,\ldots,a_q\}$, and set
\[
 \Lambda_K(p)=\sup_{q\ge p}\lambda_K(F_{p,q}).
\]

\begin{proposition}[Cost criterion]\label{prop:cost-criterion}
The condition $\Lambda_K(p)\to0$ is equivalent to the existence of a cover
of $K\setminus\{0\}$ by finite sets $F_j$ such that
$\sum_j\lambda_K(F_j)<\infty$. If these equivalent conditions hold, then
$K$ is exactly realized by a real even purely singular Rajchman measure on
$\T$.
\end{proposition}

\begin{proof}
First, $\lambda_K$ is monotone and subadditive. Monotonicity is
immediate. For subadditivity, take polynomials arbitrarily close to the
infima for finite sets $F,F'\subset K\setminus\{0\}$. Modulate one
polynomial far enough that the two Fourier supports are disjoint.
Modulation preserves its $L^1$ norm and its nonzero coefficient values,
so their sum is admissible for $F\cup F'$ and
\[
 \lambda_K(F\cup F')\le\lambda_K(F)+\lambda_K(F').
\]
Every finite set has finite cost, for example by placing its values at
distinct frequencies.

We next justify the equivalence of the two hypotheses. If $\Lambda_K(p)\to0$,
choose $1=p_0<p_1<p_2<\cdots$ so that
$\Lambda_K(p_j)\le2^{-j}$ for $j\ge1$. The initial tranche
$F_{1,p_1-1}$ has finite cost, and the remaining tranches
$F_{p_j,p_{j+1}-1}$ have summable costs. Conversely, suppose that finite
sets $F_j$ cover $K\setminus\{0\}$ and
$\sum_j\lambda_K(F_j)<\infty$. Given $\delta>0$, choose $J$ such that
$\sum_{j>J}\lambda_K(F_j)<\delta$, and choose $p$ beyond all indices
of the values in $F_1\cup\cdots\cup F_J$. Every finite tranche
$F_{p,q}$ is covered by finitely many $F_j$ with $j>J$. Monotonicity
and subadditivity give $\lambda_K(F_{p,q})<\delta$, uniformly for
$q\ge p$. Thus $\Lambda_K(p)\to0$.

Fix such a cover, indexed by $j\ge1$, and choose polynomials $R_j$
with summable errors above their infima. Then
\begin{equation}\label{eq:cost-forcing-data}
 F_j\subset\wh R_j(\Z)\subset K,
 \qquad \sum_{j\ge1}\norm{R_j}_1<\infty.
\end{equation}
Write $R_j^\vee(t)=R_j(-t)$. At stage $j$ choose a sufficiently large
integer $M_j$ and set
\begin{equation}\label{eq:cost-symmetric-forcing}
 T_j=e_{M_j}R_j+e_{-M_j}R_j^\vee.
\end{equation}
The two translated supports are disjoint and lie outside the entire
previously protected interval. Their coefficients are reflected copies
of the same nonnegative real coefficients, so $T_j$ is real and even,
all its coefficients belong to $K$, and every value in $F_j$ occurs.
Moreover,
\begin{equation}\label{eq:cost-forcing-sum}
 \norm{T_j}_1\le2\norm{R_j}_1,
 \qquad \sum_{j\ge1}\norm{T_j}_1<\infty.
\end{equation}

We interlace these forcing steps with the circle flat blocks from
Section~\ref{sec:realization}. Put $A=a_1$, start with $P_0=0$ and
$L_0=1$, and protect all frequencies in $[-2L_0,2L_0]$, in particular
frequency zero. Choose positive errors $\eta_j$ with
$\sum_j\eta_j<\infty$. Suppose that $P_{j-1}$ and $L_{j-1}$ have
been fixed. Place $T_j$ outside $[-2L_{j-1},2L_{j-1}]$ and put
$\widetilde P_j=P_{j-1}+T_j$. Choose an integer
\[
 H_j\ge\max\{2L_{j-1},\deg\widetilde P_j\}.
\]
All coefficients of $\widetilde P_j$ are in $K\subset[0,A]$. Apply
the circle flat-block construction to $\widetilde P_j$, with protected
interval $[-H_j,H_j]$ and error $\eta_j$. Choose its new coefficient
$c_j\in K\setminus\{0\}$ so small that $c_j<2^{-j}$ and the block
construction is available. For the resulting $Q_j$ set
\[
 P_j=\widetilde P_j+Q_j,
 \qquad L_j\ge\max\{L_{j-1},j,1,\deg P_j\}.
\]
All later forcing steps and flat blocks are required to avoid
$[-2L_j,2L_j]$. Thus each coefficient, including each protected zero,
is fixed after finitely many steps. All coefficients of every $P_j$
belong to $K$, and
\begin{equation}\label{eq:cost-norm-recursion}
 \norm{P_j}_1\le\norm{P_{j-1}}_1+\norm{T_j}_1+\eta_j.
\end{equation}
The sequence is bounded in $M(\T)$. Coefficient stabilization and
density of trigonometric polynomials in $C(\T)$ give a unique weak-star
limit $\mu$. Every $a_n$ occurs, by the covering property, and no
coefficient outside $K$ has been introduced. Since the coefficient at
zero stays zero, we obtain the exact equality $\wh\mu(\Z)=K$.

The Hellinger estimate and~\eqref{eq:H-continuity} give
\begin{equation}\label{eq:cost-H-recursion}
 \HH(P_j)\le\theta\HH(P_{j-1})
       +\theta\sqrt{\norm{T_j}_1}+\eta_j.
\end{equation}
The last two terms tend to zero by~\eqref{eq:cost-forcing-sum}, so the
contraction recursion implies $\HH(P_j)\to0$. The Fourier support of
$P_j$ is contained in $[-L_j,L_j]$, and all later terms avoid
$[-2L_j,2L_j]$. The multiplier identity~\eqref{eq:V-multiplier}
therefore gives the exact recovery formula
\[
 V_{L_j}*\mu=P_j.
\]
Since $L_j\to\infty$, Lemma~\ref{lem:ac-recovery} identifies the limit
in measure as the density of the absolutely continuous part of $\mu$.
On the other hand, for every $t>0$,
\[
 m_\T\{|P_j|>t\}\le t^{-1/2}\HH(P_j)\longrightarrow0.
\]
That density is consequently zero, so $\mu\perp m_\T$.

Finally, each coefficient introduced by $T_j$ is bounded in modulus by
$\norm{T_j}_1\to0$, and every coefficient of $Q_j$ is either zero or
$c_j\to0$. Each stage uses only finitely many frequencies, hence
$\mu\in M_0(\T)$. Reality and evenness pass to the weak-star limit.
\end{proof}

\begin{theorem}[Dini criterion]\label{thm:dini}
If $\sum_{n\ge1}a_n/n<\infty$, then $\Lambda_K(p)\to0$. Hence $K$ has
a real even purely singular Rajchman exact realization on $\T$.
\end{theorem}

\begin{proof}
Let $b_1\ge\cdots\ge b_N\ge0$, set $b_{N+1}=0$, and put
$D_k^+=\sum_{j=0}^{k-1}e_j$. Abel summation gives
\[
 \sum_{k=1}^Nb_ke_{k-1}
 =\sum_{k=1}^N(b_k-b_{k+1})D_k^+.
\]
Using the standard estimate $\norm{D_k^+}_1\le C\log(k+1)$, including
the terminal contribution $b_N\log(N+1)$, we obtain
\begin{align*}
 \left\|\sum_{k=1}^Nb_ke_{k-1}\right\|_1
 &\le C\sum_{k=1}^N(b_k-b_{k+1})\log(k+1)\\
 &=C\left(b_1\log2+\sum_{k=2}^N
 b_k\bigl(\log(k+1)-\log k\bigr)\right)\\
 &\le C'\sum_{k=1}^N\frac{b_k}{k}.
\end{align*}
Apply this to $b_k=a_{p+k-1}$. Then
\[
 \lambda_K(F_{p,q})\le C\sum_{k\ge1}\frac{a_{p+k-1}}{k}.
\]
The right-hand side tends to zero by dominated convergence, dominated
by $a_k/k$.
\end{proof}

\begin{theorem}[Palettes with discretely convex enumerations]\label{thm:convex}
If $a_n-2a_{n+1}+a_{n+2}\ge0$ for all $n$, then
\[
 \lambda_K(F_{p,q})=a_p,\qquad \Lambda_K(p)=a_p.
\]
Consequently, every positive null palette whose decreasing enumeration is
discretely convex has a real even purely singular Rajchman exact
realization on $\T$, and its norm may be made
arbitrarily close to $a_1$.
\end{theorem}

\begin{proof}
For the tail $b_n=a_{p+n}$, $n\ge0$, let $d_n=b_n-b_{n+1}$ and
$w_j=(j+1)(d_j-d_{j+1})$. Discrete convexity makes $(d_n)$ nonnegative
and decreasing, and $\sum_{n\ge0}d_n=b_0$. In particular,
$(j+1)d_{j+1}\to0$. Summation by parts gives
\[
 w_j\ge0,\qquad
 \sum_{j=0}^{N}w_j=b_0-b_{N+1}-(N+1)d_{N+1}
 \longrightarrow b_0.
\]
The Fej\'er kernels used to define $V_N$ above give a
nonnegative $L^1$ function
\[
 f_b=\sum_{j\ge0}w_j\mathsf F_{j+1},
 \qquad \norm{f_b}_1=b_0.
\]
This is the classical convex-sequence representation by Fej\'er kernels;
see \cite[Theorem~I.4.1, p.~23]{Katznelson2004}. The coefficient
quantization below is the additional step needed for exact palettes.
For $k\ge0$, their Fourier coefficients and another telescoping give
\[
 \wh f_b(k)=\wh f_b(-k)
 =\sum_{j\ge k}(j+1-k)(d_j-d_{j+1})
 =\sum_{j\ge k}d_j=b_k.
\]

Let $g_N=\mathsf F_{N+1}*f_b$. For $|n|\le N$ its coefficient is
\[
 \wh g_N(n)=\left(1-\frac{|n|}{N+1}\right)b_{|n|}>0.
\]
Replace this coefficient by the smallest value in the tail palette
$\{b_0,b_1,\ldots\}$ which is not smaller than it, and set all
coefficients outside $[-N,N]$ equal to zero. Denote the resulting real
even polynomial by $Q_N$. If the chosen value is $b_m$, then
$m\ge|n|$ because $\wh g_N(n)\le b_{|n|}$, and
\[
 0\le\wh Q_N(n)-\wh g_N(n)
   <b_m-b_{m+1}=d_m\le d_{|n|}.
\]
For fixed $n$, the inequality
$\wh g_N(n)>b_{|n|+1}$ holds for all sufficiently large $N$.
Consequently $\wh Q_N(n)=b_{|n|}$ eventually, and the rounding error
then equals $|n|b_{|n|}/(N+1)$, which tends to zero. Since
$\sum_{n\in\Z}d_{|n|}^2<\infty$, Parseval and dominated convergence
give
\[
 \norm{Q_N-g_N}_2\longrightarrow0,
 \qquad \norm{Q_N-g_N}_1\longrightarrow0.
\]
Also $g_N\to f_b$ in $L^1$, so $\norm{Q_N}_1\to a_p$.
For each fixed $q\ge p$, all values $a_p,\ldots,a_q$ occur among the
coefficients of $Q_N$ for all sufficiently large $N$. Thus
$\lambda_K(F_{p,q})\le a_p$. The reverse inequality follows from
$|\wh P(n)|\le\norm P_1$, proving both cost identities.
Proposition~\ref{prop:cost-criterion} gives the asserted exact singular
Rajchman realization.

For the near-optimal norm, fix $\delta>0$. Choose
$2\le p_1<p_2<\cdots$ and positive numbers $\zeta_j$ such that
\[
 2\sum_{j\ge1}(a_{p_j}+\zeta_j)<\frac\delta2.
\]
Use the preceding quantization with $p=1$ to choose a real even
polynomial $P_0$ realizing $a_1,\ldots,a_{p_1-1}$, with all
coefficients in $K$, $\wh P_0(0)=a_1$, and
$\norm{P_0}_1<a_1+\delta/4$. Choose a nonzero integer $\xi$ outside
its Fourier support and an initial protected interval containing both
that support and $\xi$; the zero coefficient at $\xi$ will remain
fixed. For each $j\ge1$, the cost identity supplies a polynomial $R_j$
realizing $F_{p_j,p_{j+1}-1}$ with
$\norm{R_j}_1<a_{p_j}+\zeta_j$. Its separated reflected modulation
$T_j$ from~\eqref{eq:cost-symmetric-forcing} satisfies
\[
 \sum_{j\ge1}\norm{T_j}_1<\frac\delta2.
\]
Run the same forcing and flat-block iteration, now from this $P_0$,
with $\sum_j\eta_j<\delta/4$. The preceding proof of exact range,
pure singularity and the Rajchman property is unchanged; the protected
zero is at $\xi$ instead of at zero. Equation~\eqref{eq:cost-norm-recursion}
gives
\[
 a_1\le\norm\mu
 \le\norm{P_0}_1+\sum_j\norm{T_j}_1+\sum_j\eta_j
 <a_1+\delta.\qedhere
\]
\end{proof}

\begin{theorem}[Sharp summability threshold]\label{thm:sharp-lr}
For $r>2$, there exists a nonzero purely singular Rajchman measure
$\mu$ on $\T$ with $\wh\mu(\Z)=K$ and $\wh\mu\in\ell^r(\Z)$ if and
only if $(a_n)\in\ell^r(\N)$. Moreover, there exists a single nonzero
purely singular Rajchman measure $\mu$ with $\wh\mu(\Z)=K$ and
$\wh\mu\in\bigcap_{r>2}\ell^r(\Z)$ if and only if
$(a_n)\in\bigcap_{r>2}\ell^r(\N)$. For $0<r\le2$ no nonzero purely
singular measure can have an $\ell^r$ transform.
\end{theorem}

\begin{proof}
Necessity follows by choosing pairwise distinct frequencies $n_k$ at
which the distinct values $a_k$ occur. Then
\[
 \sum_{k\ge1}a_k^r\le\sum_{n\in\Z}|\wh\mu(n)|^r.
\]

For sufficiency, fix $r>2$ and assume $(a_n)\in\ell^r$. H\"older's
inequality with conjugate exponent $r'=r/(r-1)>1$ gives
\[
 \sum_{n\ge1}\frac{a_n}{n}
 \le\left(\sum_{n\ge1}a_n^r\right)^{1/r}
     \left(\sum_{n\ge1}n^{-r'}\right)^{1/r'}<\infty.
\]
We use the explicit one-sided polynomials from the proof of
Theorem~\ref{thm:dini}, rather than arbitrary nearly optimal forcing
polynomials:
\begin{equation}\label{eq:lr-explicit-forcing}
 R_{p,q}(t)=\sum_{k=1}^{q-p+1}a_{p+k-1}e_{k-1}(t).
\end{equation}
Write $B(p)=C\sum_{k\ge1}a_{p+k-1}/k$. Then $B(p)\to0$ and
$\norm{R_{p,q}}_1\le B(p)$ uniformly in $q\ge p$. Choose
$1=p_0<p_1<p_2<\cdots$ so that $B(p_j)\le2^{-j}$ for $j\ge1$.
At stage $j\ge1$ use $R_{p_{j-1},p_j-1}$ and its reflected copy,
modulated into disjoint intervals beyond all protected frequencies.
The initial forcing polynomial has finite norm, and the subsequent
forcing norms are summable.

The tranches partition the prescribed values, and the one-sided
polynomial~\eqref{eq:lr-explicit-forcing} contains each value in its
tranche exactly once. Hence, after reflection, every $a_n$ occurs
exactly twice in the forcing part. If $T_j$ denotes the resulting
forcing polynomial, the disjointness of all these spectra gives
\begin{equation}\label{eq:lr-forcing-budget}
 \sum_{j\ge1}\sum_{n\in\supp\wh T_j}|\wh T_j(n)|^r
 =2\sum_{n\ge1}a_n^r<\infty.
\end{equation}

Interlace the forcing steps with the flat blocks exactly as in
Proposition~\ref{prop:cost-criterion}, using the uniform coefficient
bound $A=a_1$. After $\widetilde P_j=P_{j-1}+T_j$ has been fixed,
let $D_j$ be the size of the old frequency interval used in
Lemma~\ref{lem:product}. Choose the new coefficient
$c_j\in K\setminus\{0\}$ small enough that the block construction
applies and
\[
 0<c_j<\min(1,2^{-j}),\qquad D_jc_j^{1/j}\le2^{-j}.
\]
The support estimate~\eqref{eq:block-count} gives
\begin{equation}\label{eq:lr-flat-budget}
 \sum_{n\in\supp\wh Q_j}|\wh Q_j(n)|^r
 \le\frac{A^2}{4}D_jc_j^{r-2}
 \le\frac{A^2}{4}2^{-j}
\end{equation}
for all sufficiently large $j$, since then $1/j\le r-2$ and
$0<c_j<1$. The finitely many remaining blocks have finite supports.
All forcing and flat-block spectra are pairwise disjoint, so
\eqref{eq:lr-forcing-budget} and~\eqref{eq:lr-flat-budget} prove
$\wh\mu\in\ell^r(\Z)$. Exact range, the Rajchman property and pure
singularity follow from the same protected-interval iteration and
exact recovery used in Proposition~\ref{prop:cost-criterion}.

If $(a_n)\in\bigcap_{r>2}\ell^r$, make the same choices of tranches
and coefficients $c_j$. For each fixed $r>2$, the forcing sum is finite
and~\eqref{eq:lr-flat-budget} applies after finitely many stages.
Thus one measure works simultaneously for every $r>2$.

Finally, if $0<r\le2$ and $\wh\mu\in\ell^r$, then
$\wh\mu\in\ell^2$. Plancherel gives an $L^2$ function with the same
Fourier coefficients as $\mu$, and uniqueness of Fourier transforms
identifies $\mu$ with that absolutely continuous measure. This excludes
every nonzero purely singular measure.
\end{proof}

\section{Applications and problems from the literature}\label{sec:applications}

The exact-palette results separate three phenomena which are often close
in the Wiener--Pitt literature: containment of the Fourier--Stieltjes
range in a thin set, exact realization of every point of that set, and
regularity of convolution powers. We now record the applications of the
block method which answer previously stated questions or sharpen their
interpretation.
\subsection{Complex ranges and open range conditions}\label{sec:complex}

The exact block structure retained in~\eqref{eq:flat-data} and
\eqref{eq:weakstar} allows a small $L^1$ correction of phases.
This step only needs a summable deviation of the selected phases
from a single limiting direction.

\begin{proof}[Proof of Theorem~\ref{thm:complex}]
Compactness of the unit circle provides a subsequence of
$z_n/|z_n|$ converging to some $u$ with $|u|=1$.
Pass to a further subsequence and relabel its terms so that
\begin{equation}\label{eq:phase-selection}
 z_n=a_nu_n,\qquad
 a_1>a_2>\cdots>0,\qquad a_n\to0,\qquad
 |u_n|=1,
\end{equation}
and
\begin{equation}\label{eq:phase-budget}
 A:=a_1<R/3,
 \qquad \sum_{n\ge1}|u_n-u|\le\frac14.
\end{equation}
These requirements can be imposed simultaneously. At each selection
take an index sufficiently far along the convergent subsequence to
decrease the modulus strictly and make the phase error smaller than
a chosen summable bound. Starting sufficiently far out also ensures
$a_1<R/3$.

Apply Theorem~\ref{thm:main} and its construction to this positive
sequence with $\eps=A$. It gives a nonzero singular measure $\mu$,
polynomials $P_j=A+\sum_{i=1}^jQ_i$, and pairwise disjoint finite
sets $E_j$ such that
\begin{equation}\label{eq:complex-blocks}
 \norm\mu\le2A,\qquad \norm{Q_j}_1\le4A,
 \qquad
 \wh Q_j=a_{n_j}\one_{E_j},
 \qquad
 \wh\mu=A\one_{\{0\}}+\sum_{j\ge1}a_{n_j}\one_{E_j},
\end{equation}
where the indices $n_j\ge2$ are strictly increasing.
As observed after~\eqref{eq:increment-norm}, these structural
assertions hold on the original group after Haar lifting as well.

The series
\begin{equation}\label{eq:phase-correction}
 f=(u_1-u)A+\sum_{j\ge1}(u_{n_j}-u)Q_j
\end{equation}
converges absolutely in $L^1(G)$, because
\begin{equation}\label{eq:phase-norm}
 \norm f_1
 \le A|u_1-u|+4A\sum_{j\ge1}|u_{n_j}-u|
 \le4A\sum_{n\ge1}|u_n-u|\le A.
\end{equation}
Define $\nu=u\mu+fm_G$.
At $0$ its Fourier coefficient is $z_1$; on $E_j$ it is $z_{n_j}$;
elsewhere it is zero. Thus its Fourier--Stieltjes range is contained in the set consisting of
zero and the selected values, and hence in the original prescribed set.
Moreover,
\begin{equation}\label{eq:modulus-preservation}
 |\wh\nu(\gamma)|=\wh\mu(\gamma)
 \qquad(\gamma\in\Gamma).
\end{equation}
It follows that $\nu$ is Rajchman and its transform belongs to
every $\ell^r$, $r>2$. The argument following~\eqref{eq:ell-r} gives
the stated regularity of the second and third convolution powers, and
Lemmas~\ref{lem:Rajchman-band} and~\ref{lem:natural-spectrum} give strong
continuity and natural spectrum.
The singular part of $\nu$ is exactly $u\mu$, which is nonzero.
Finally, $\norm\nu\le3A<R$ by~\eqref{eq:phase-norm}.
\end{proof}

We now replace each remaining zero Fourier coefficient by a nonzero value
in the prescribed open set, using an absolutely convergent Fourier series.
Only this step requires countability of the dual.

\begin{proof}[Proof of Theorem~\ref{thm:open}]
Since $U$ is open and $0\in\overline U$, choose nonzero
$z_n\in U$ with $z_n\to0$.
Theorem~\ref{thm:complex}, with norm bound $R/2$, gives
$\rho\in M_0(G)\setminus L^1(G)$ with Fourier--Stieltjes range contained in
$\{0,z_1,z_2,\ldots\}$ and all the asserted summability properties.
If $0\in U$, take $\nu=\rho$.

Suppose $0\notin U$. Compact metrizability implies that $\Gamma$
is countable. For completeness, its characters form an orthonormal
set in the separable Hilbert space $L^2(G)$, so there can only be
countably many of them.
Enumerate the zero set of $\wh\rho$ as
$\{\gamma_k:k\in I\}$, where $I$ is a finite initial segment of
$\N$ or $I=\N$; the empty set is permitted.
Choose $w_k\in U$ with
\[
 |w_k|<R2^{-k-2}.
\]
The series
\[
 g(x)=\sum_{k\in I}w_k\gamma_k(x)
\]
converges absolutely and uniformly. Its transform equals $w_k$ at
$\gamma_k$ and zero elsewhere, by character orthogonality.
In particular,
\[
 \norm g_1\le\norm g_\infty\le\sum_{k\in I}|w_k|<R/2,
 \qquad \wh g\in\ell^1(\Gamma).
\]
Set $\nu=\rho+gm_G$. All old nonzero Fourier coefficients remain
unchanged and belong to $U$, and each old zero is replaced by the
corresponding $w_k\in U$. Thus $\wh\nu(\Gamma)\subset U$.
The singular parts of $\nu$ and $\rho$ coincide and are nonzero,
and $\norm\nu<R$.
The $\ell^1$ correction preserves the Rajchman condition and every
$\ell^r$ property for $r>2$. The convolution-power and spectral
conclusions follow as before.
\end{proof}

\begin{remark}\label{rem:countability}
The metrizability assumption in Theorem~\ref{thm:open} is needed only
when $0\notin U$, and in that case it has a concrete role. A function
in $c_0(\Gamma)$ is supported on
a countable set, since its support is the union of the finite sets
where its modulus is at least $1/n$. Hence a Rajchman transform on
an uncountable dual must vanish somewhere, and cannot have its
entire range in such a $U$. Theorem~\ref{thm:complex}, which permits
the value zero, remains valid without metrizability.
\end{remark}

\begin{remark}
The measures furnished by Theorem~\ref{thm:open} are continuous and
not absolutely continuous, although their convolution squares belong to
$L^p(G)$ for every $1\le p<\infty$ and their third convolution powers
have continuous densities with absolutely convergent Fourier series. Thus no
open set accumulating at zero can force absolute
continuity of every continuous measure on the circle with Fourier
range in that set. This is the absolute-continuity strengthening
asked about in \cite[Section~7, item~2]{OW2014}.
The conclusion concerns the measure itself. It is compatible with
the natural-spectrum conclusions of \cite{OSW2026,OW2014}, since
all the measures constructed here have natural spectra.

The condition $0\in\overline U$ excludes only vacuous cases of the
question on the circle. For a continuous complex measure $\rho$,
Wiener's mean-square formula gives
\[
 \lim_{N\to\infty}\frac1{2N+1}
     \sum_{|n|\le N}|\wh\rho(n)|^2
 =\sum_{x\in\T}|\rho(\{x\})|^2=0.
\]
One obtains the identity by expressing the mean as a double integral
against $\rho\otimes\overline\rho$: the normalized Dirichlet
kernels are bounded in modulus by one and converge pointwise to
the indicator of the diagonal. Dominated convergence then gives
the sum of squared atom masses. If $0\notin\overline U$, some
$\delta>0$ satisfies $|z|\ge\delta$ for every $z\in U$, which is
incompatible with this formula and $\wh\rho(\Z)\subset U$.
\end{remark}

\subsection{A phase-valued construction for uniform convergence}
\label{sec:phase-application}

We now apply the block method on the circle in a different direction.
The Boolean selector from Lemma~\ref{lem:selector} is replaced by a
selector with values among the fourth roots of unity. Its Walsh coefficient
indexed by the full set $[d]$ reproduces the old polynomial, while its
pointwise square is the parity character
$\chi(s)=s_1\cdots s_d$. After the path substitution, the latter becomes
the Rudin--Shapiro sign, so squaring the Fourier coefficients of a block
recovers those signs exactly.
This additional cancellation gives uniform convergence of the Fourier
series of the convolution square.

\subsubsection{Fourier partial sums and two obstructions}

For $f\in C(\T)$, allowing the value $+\infty$, put
\begin{equation}\label{eq:U-norm}
 S_Nf(t)=\sum_{|n|\le N}\wh f(n)e_n(t),
 \qquad
 \norm f_U=\sup_{N\ge0}\norm{S_Nf}_\infty.
\end{equation}
Thus $S_Nf$ is the symmetric Fourier partial sum of order $N$. We write
\[
 U(\T)=\{f\in C(\T):S_Nf\to f\text{ uniformly}\}
\]
for the space of continuous functions with uniformly convergent Fourier
series; $\norm{\cdot}_U$ is finite on this space. We also write
\[
 A(\T)=\left\{f\in C(\T):
             \sum_{n\in\Z}|\wh f(n)|<\infty\right\}
\]
for the Wiener algebra.

By contrast, the displayed definition of the $u$-norm in
\cite[p.~10]{GuptaHare} takes, for a compact group, the supremum
\begin{equation}\label{eq:GH-u-norm}
 \norm f_u=
 \sup_{\mathcal F\Subset\wh G}
 \left\|
   \sum_{\pi\in\mathcal F}d_\pi
   \operatorname{Tr}\bigl(\wh f(\pi)\pi(\,\cdot\,)\bigr)
 \right\|_\infty,
\end{equation}
over all finite subsets of the unitary dual. The proof on
pp.~12--13 instead uses initial segments of a fixed enumeration of the
dual \cite[pp.~12--13]{GuptaHare}. The distinction is
essential already on the circle.

\begin{proposition}\label{prop:unconditional-obstruction}
Let $\rho\in M(\T)$. If
\begin{equation}\label{eq:finite-set-bounded}
 \sup_{F\Subset\Z}
 \left\|\sum_{n\in F}\wh\rho(n)^2e_n\right\|_\infty<\infty,
\end{equation}
then $\rho$ has an $L^2$ density. In particular, no nonzero singular
measure can satisfy the literal finite-set condition
\eqref{eq:finite-set-bounded}.
\end{proposition}

\begin{proof}
Evaluate at $t=0$. The finite subsums of the scalar family
$(\wh\rho(n)^2)_{n\in\Z}$ are uniformly bounded. Cover the complex
plane by eight closed sectors of angle $\pi/4$, assigning
boundary points to one of the adjacent sectors. In each sector choose
a unimodular number $\omega$ such that
\[
 \RePart(\omega z)\ge \cos(\pi/8)|z|
\]
for every $z$ in that sector. Applying the bound to arbitrary finite
subsets contained in one sector shows that the sum of the moduli in
that sector is finite. Summing over the eight sectors gives
\[
 \sum_{n\in\Z}|\wh\rho(n)|^2<\infty.
\]
By Plancherel there is $g\in L^2(\T)$ with
$\wh g=\wh\rho$, and uniqueness of Fourier--Stieltjes transforms gives
$\rho=g\,m_{\T}$.
\end{proof}

A nonzero purely singular even measure whose convolution square has a
uniformly convergent Fourier series cannot be real.

\begin{proposition}\label{prop:real-even-obstruction}
Let $\rho\in M(\T)$ be real and even. If
\[
 \rho*\rho=fm_{\T}
 \qquad\text{for some }f\in L^\infty(\T),
\]
then $\rho\in L^2(\T)$ and
\[
 \sum_{n\in\Z}|\wh\rho(n)|^2\le\norm f_\infty.
\]
Hence no nonzero real even purely singular measure can have a
convolution-square density in $L^\infty(\T)$. In particular, it cannot have
a uniformly convergent convolution-square Fourier series.
\end{proposition}

\begin{proof}
Reality and evenness give
$\wh\rho(-n)=\wh\rho(n)\in\mathbb R$, and therefore
\[
 \wh f(n)=\wh\rho(n)^2=|\wh\rho(n)|^2.
\]
Let $\mathsf F_N$ be the normalized Fej\'er kernel with Fourier
coefficients $1-|n|/N$ for $|n|<N$. Positivity and normalization give
\begin{align*}
 0&\le\sum_{|n|<N}\left(1-\frac{|n|}{N}\right)
       |\wh\rho(n)|^2\\
  &=(\mathsf F_N*f)(0)\le\norm f_\infty.
\end{align*}
The summands increase pointwise to $|\wh\rho(n)|^2$. Monotone
convergence yields the displayed $\ell^2$ estimate. Plancherel gives an
$L^2$ function with Fourier transform $\wh\rho$, and uniqueness of
Fourier--Stieltjes transforms identifies it with $\rho$.
\end{proof}

\begin{remark}\label{rem:el-abdalaoui}
The published abstract of \cite{ElAbdalaoui2023} announces a symmetric
singular probability measure whose convolution square has a continuous
nearly flat density. If ``symmetric'' is understood as reflection
invariance, the announced conclusion conflicts with
Proposition~\ref{prop:real-even-obstruction}, since a continuous density is
bounded. The abstract alone does not provide enough detail to resolve this
discrepancy, and we do not use the claim as an antecedent. The phase-valued
construction below avoids the obstruction because its even Fourier
coefficients need not be real.
\end{remark}

Thus Theorem~\ref{thm:uniform-square} is necessarily a genuinely
complex-valued statement. The construction below is even under
reflection, but it is not real.

\subsubsection{The phase-valued selector}

Put
\[
 \Omega_4=\{1,-1,i,-i\},
 \qquad
 \chi(s)=s_1\cdots s_d
 \quad(s\in\{-1,1\}^d),
\]
and, for $d\ge2$, let
\begin{equation}\label{eq:phase-grid}
 \mathcal G_d=
 \bigl(2^{1-d}\Z\cap[-1/2,1/2]\bigr)
 +i\bigl(2^{1-d}\Z\cap[-1/2,1/2]\bigr).
\end{equation}

\begin{lemma}[Phase selector]\label{lem:phase-selector}
For every $b\in\mathcal G_d$ there is a function
$F_b:\{-1,1\}^d\to\Omega_4$ such that
\begin{equation}\label{eq:phase-selector-properties}
 F_b(s)^2=\chi(s),
 \qquad
 \wh F_b([d])=b.
\end{equation}
Moreover,
\begin{equation}\label{eq:phase-Walsh-l1}
 \sum_{S\subset[d]}|\wh F_b(S)|\le2^{d/2}.
\end{equation}
Every $z\in\C$ with $|z|\le1/2$ is within distance
$\sqrt2\,2^{-d}$ of some point of $\mathcal G_d$.
\end{lemma}

\begin{proof}
Write $b=u+iv$ and split the cube into
$C_+=\{s:\chi(s)=1\}$ and $C_-=\{s:\chi(s)=-1\}$.
Both classes have $2^{d-1}$ elements. On $C_+$ choose values
$F_b(s)\in\{1,-1\}$ so that
\[
 \sum_{s\in C_+}F_b(s)=2^du.
\]
On $C_-$ put $F_b(s)=i\epsilon_s$, where
$\epsilon_s\in\{-1,1\}$ are chosen so that
\[
 \sum_{s\in C_-}\epsilon_s=-2^dv.
\]
These choices are possible because sums of $2^{d-1}$ signs range from
$-2^{d-1}$ to $2^{d-1}$ in steps of two, exactly as prescribed by
\eqref{eq:phase-grid}. We then have $F_b^2=1$ on $C_+$ and
$F_b^2=-1$ on $C_-$, proving the first identity. Also,
\[
 \wh F_b([d])
 =2^{-d}\left(
   \sum_{C_+}F_b(s)-i\sum_{C_-}\epsilon_s
  \right)=u+iv=b.
\]
Since $|F_b|=1$, Parseval on the Boolean cube and
Cauchy--Schwarz give
\[
 \sum_{S\subset[d]}|\wh F_b(S)|
 \le2^{d/2}
    \left(\sum_{S\subset[d]}|\wh F_b(S)|^2\right)^{1/2}
 =2^{d/2}.
\]
Finally, round the real and imaginary parts of $z$ independently to
the grid of mesh $2^{1-d}$. Each coordinate changes by at most
$2^{-d}$.
\end{proof}

Using the paired colouring~\eqref{eq:colours}, define
\begin{equation}\label{eq:phase-polynomial}
 \mathcal B_b(z)=
 \sum_{x\in\{0,1\}^m}
 F_b\bigl(\sigma_1(x),\ldots,\sigma_d(x)\bigr)z^x.
\end{equation}
The same pointwise Walsh expansion as in Lemma~\ref{lem:selector}
gives
\begin{equation}\label{eq:phase-expansion}
 \mathcal B_b
 =bR_m+\sum_{S\subsetneq[d]}\wh F_b(S)R_{m,S}.
\end{equation}
Consequently, Lemma~\ref{lem:proper} and
\eqref{eq:phase-Walsh-l1} imply
\begin{equation}\label{eq:phase-selector-error}
 2^{-m/2}\norm{\mathcal B_b-bR_m}_1
 \le2^{d/2}q^{r(m,d)}.
\end{equation}
The feature not present in the Boolean selector is the exact identity
\begin{equation}\label{eq:phase-square-RS}
 F_b\bigl(\sigma_1(x),\ldots,\sigma_d(x)\bigr)^2
 =\prod_{j=1}^d\sigma_j(x)
 =(-1)^{x_1x_2+\cdots+x_{m-1}x_m}.
\end{equation}
Thus squaring the coefficients of $\mathcal B_b$ removes the parameter
$b$ and leaves precisely the Rudin--Shapiro signs.

\subsubsection{A maximal Rudin--Shapiro estimate}

For $x\in\{0,1\}^m$ put
\[
 \nu(x)=x_1+2x_2+\cdots+2^{m-1}x_m,
 \qquad
 \epsilon_x=(-1)^{x_1x_2+\cdots+x_{m-1}x_m}.
\]

\begin{lemma}\label{lem:RS-prefix}
For every $m\ge1$,
\begin{equation}\label{eq:RS-prefix}
 \max_{r\in\{0,\ldots,2^m\}}
 \left\|
  \sum_{\nu(x)<r}\epsilon_xz^x
 \right\|_\infty
 \le(2+\sqrt2)2^{m/2}.
\end{equation}
\end{lemma}

\begin{proof}
Let $R_0=T_0=1$ and use the classical recursion
\[
 R_m=R_{m-1}+z_mT_{m-1},
 \qquad
 T_m=R_{m-1}-z_mT_{m-1}.
\]
The polynomial $R_m$ is exactly the one in~\eqref{eq:RS}, with its
monomials ordered by $\nu(x)$. Let $R_m^{(r)}$ and $T_m^{(r)}$ denote
the first $r$ monomials in this binary order, and put
\[
 M_m=\max_{0\le r\le2^m}
 \max\bigl(\norm{R_m^{(r)}}_\infty,
           \norm{T_m^{(r)}}_\infty\bigr).
\]
For a prefix lying in the first half, the relevant polynomial is a
prefix from level $m-1$. For $0\le s\le2^{m-1}$,
\[
 R_m^{(2^{m-1}+s)}=R_{m-1}+z_mT_{m-1}^{(s)},
 \qquad
 T_m^{(2^{m-1}+s)}=R_{m-1}-z_mT_{m-1}^{(s)}.
\]
The complementary identity gives
$\norm{R_{m-1}}_\infty\le2^{m/2}$. Hence
\[
 M_m\le2^{m/2}+M_{m-1}.
\]
After setting $u_m=2^{-m/2}M_m$, we obtain
$u_m\le1+2^{-1/2}u_{m-1}$. Since $u_0=1$, induction gives
$u_m\le(1-2^{-1/2})^{-1}=2+\sqrt2$. The estimate for the prefixes of
$R_m$ is~\eqref{eq:RS-prefix}.
\end{proof}

\subsubsection{A complex even block on the circle}

We now record the phase version of the circle block. Only the new
points are proved in detail; the $L^1$ and Hellinger estimates follow
the argument of Lemma~\ref{lem:product} and the circle realization in
Section~\ref{sec:realization}.

\begin{lemma}\label{lem:phase-block}
Let
\[
 P(t)=\sum_{|k|\le F}p_ke_k(t),
 \qquad p_{-k}=p_k,
 \qquad |p_k|\le A,
\]
where $A>0$ and $F\ge0$ is an integer, and put $D=2F+1$. Let
$H\ge F$ be an integer and let $\eta>0$.
For every sufficiently small $a>0$ there is an even complex
trigonometric polynomial $Q$ such that:
\begin{enumerate}[label=\textnormal{(\roman*)},leftmargin=2.2em]
\item $\supp\wh Q\cap[-H,H]=\varnothing$ and every nonzero Fourier
coefficient of $Q$ belongs to $a\Omega_4$;
\item
\begin{equation}\label{eq:phase-block-H}
 \norm{P+Q}_1\le\norm P_1+\eta,
 \qquad
 \HH(P+Q)\le\theta\HH(P)+\eta;
\end{equation}
\item
\begin{equation}\label{eq:phase-block-count}
 \#\supp\wh Q\le\frac{DA^2}{4a^2};
\end{equation}
\item
\begin{equation}\label{eq:phase-block-U}
 \norm{Q*Q}_U\le3ADa.
\end{equation}
\end{enumerate}
Here $\theta=1-1/1024$ is the same constant as in
Lemma~\ref{lem:product}.
\end{lemma}

\begin{proof}
Choose $d$ and, for each $|k|\le F$, choose
$b_k=b_{-k}\in\mathcal G_d$ so that
\begin{equation}\label{eq:phase-bk}
 \left|b_k-\frac{p_k}{2A}\right|
 \le\sqrt2\,2^{-d}.
\end{equation}
Put
\[
 S(t,z)=\sum_{|k|\le F}e_k(t)\mathcal B_{b_k}(z),
 \qquad
 S^\vee(t,z)=S(-t,z^{-1}),
\]
and, on $\T\times\T^m\times\T$, define
\begin{equation}\label{eq:phase-product-block}
 \mathcal Q(t,z,w)=a\bigl(wS(t,z)+w^{-1}S^\vee(t,z)\bigr).
\end{equation}
The second term is a reflected, not a conjugated, copy of the first.
Thus the coefficients at opposite frequencies are equal and
$\mathcal Q$ is even. Every nonzero coefficient belongs to
$a\Omega_4$.

Let $M=2^m$, $s=a\sqrt M$, and choose
\begin{equation}\label{eq:phase-m-choice}
 m=\left\lfloor2\log_2\frac{A}{4a}\right\rfloor,
 \qquad
 \beta=\frac{s}{A}.
\end{equation}
For all sufficiently small $a$ this is a positive integer and
\begin{equation}\label{eq:phase-beta}
 \frac1{4\sqrt2}<\beta\le\frac14.
\end{equation}
Write
\[
 A_d(t)=\sum_{|k|\le F}b_ke_k(t),
 \qquad
 Z_m=M^{-1/2}R_m,
 \qquad
 X=\RePart(wZ_m).
\]
Since $A_d$ is even, the full Walsh components in
\eqref{eq:phase-product-block} sum to $2sA_dX$. By
\eqref{eq:phase-bk}, \eqref{eq:phase-selector-error}, and the triangle
inequality,
\begin{equation}\label{eq:phase-product-error}
 \norm{\mathcal Q-\beta PX}_1
 \le C AD\left(2^{-d}+2^{d/2}q^{r(m,d)}\right)
\end{equation}
with an absolute constant $C$.

The random variable $X$ has the same properties as in
\eqref{eq:X-moments} for $\kappa=2$:
\[
 \E X=0,\qquad \E X^2=\frac12,
 \qquad |X|\le\sqrt2.
\]
In particular, $1+\beta X>0$. Although $P$ is now complex,
\[
 |P(1+\beta X)|=|P|(1+\beta X).
\]
Consequently the ideal norm identity and the contraction estimate
\eqref{eq:ideal} remain valid. Denote the right-hand side of
\eqref{eq:phase-product-error} by $\delta$. Choose $d$ first and then $a$
small enough that
\[
 \delta\le\min(\eta/2,\eta^2/4).
\]
In the finite averaging step below choose a common error
$\tau\le\eta/2$ for the integrals of $|P+\mathcal Q|$ and
$\sqrt{|P+\mathcal Q|}$. Then
$\delta+\tau\le\eta$ and $\sqrt\delta+\tau\le\eta$, so
\eqref{eq:H-continuity} proves~\eqref{eq:phase-block-H}.

We realize the product block on the circle exactly as in
Section~\ref{sec:realization}, with this averaging tolerance $\tau$, using
the separation conditions already imposed in
\eqref{eq:circle-separation} (with $F_0=F$):
\begin{equation}\label{eq:phase-frequency-order}
 N_j>2F+\sum_{\ell<j}N_\ell\quad(1\le j\le m),
 \qquad
 L>H+F+\sum_{j=1}^mN_j.
\end{equation}
Besides separating the support, these conditions ensure the required
binary ordering of the frequencies. Together with the finite-averaging
lacunarity conditions, they are all lower bounds and can therefore be
imposed simultaneously. After the
substitution $z_j=e_{N_j}$ and $w=e_L$, the positive frequencies are
\begin{equation}\label{eq:phase-frequencies}
 L+k+\sum_{j=1}^m x_jN_j,
 \qquad |k|\le F,
 \qquad x\in\{0,1\}^m.
\end{equation}
They are distinct, exceed $H$, and are ordered first by the binary
number $\nu(x)$ and then by $k$. The reflected copy gives their
negatives. Thus the coefficient at the frequency in
\eqref{eq:phase-frequencies} is
\[
 aF_{b_k}\bigl(\sigma_1(x),\ldots,\sigma_d(x)\bigr).
\]
The support and coefficient assertions follow, and
$\#\supp\wh Q\le2DM\le DA^2/(4a^2)$ by
\eqref{eq:phase-m-choice}.

By~\eqref{eq:phase-square-RS}, the coefficient of $Q*Q$ at the same
positive frequency is
\begin{equation}\label{eq:phase-square-coefficient}
 a^2\epsilon_x,
\end{equation}
independently of $k$. A symmetric partial sum contains some number of
complete packets indexed by $x$, at most one incomplete packet, and
the reflected negative frequencies. A complete packet contributes the
Dirichlet factor $\sum_{|k|\le F}e_k$, whose uniform norm is at most
$D$. Lemma~\ref{lem:RS-prefix} therefore gives
\begin{equation}\label{eq:phase-U-estimate}
 \norm{Q*Q}_U
 \le2Da^2\bigl((2+\sqrt2)\sqrt M+1\bigr).
\end{equation}
Since $a\sqrt M\le A/4$ and, for sufficiently small $a$,
$a\le A/4$, the right-hand side is at most $3ADa$.
\end{proof}

\subsubsection{Iteration and proof of the application}

\begin{proof}[Proof of Theorem~\ref{thm:uniform-square}]
Pass to a strictly decreasing subsequence of the prescribed null
sequence whenever necessary. Choose positive $\eta_j$ with
$\sum_j\eta_j\le\eps$, and start with $P_0=A$. We construct complex
even polynomials
\[
 P_j=P_{j-1}+Q_j
\]
using Lemma~\ref{lem:phase-block}. Set $L_0=1$. As in the
circle part of Section~\ref{sec:iteration}, before stage $j$ choose an
integer $F_j$ such that
\[
 \supp\wh P_{j-1}\subset[-F_j,F_j],
 \qquad F_j\le L_{j-1},
\]
put $D_j=2F_j+1$, and apply the block lemma with protected radius
$H=2L_{j-1}$.

Choose the next subsequence element $c_j$ so far out that the block
lemma is available with $a=c_j$ and error $\eta_j$, and also
\begin{equation}\label{eq:phase-c-choice}
 0<c_j<\min(A,1,2^{-j},c_{j-1}/2),
 \qquad
 D_jc_j^{1/j}\le2^{-j},
\end{equation}
where $c_0=A$. After defining $P_j=P_{j-1}+Q_j$, choose
\[
 L_j\ge\max(L_{j-1},j,1,\deg P_j).
\]
The support of $Q_j$ misses $[-2L_{j-1},2L_{j-1}]$, so the intervals
$[-2L_j,2L_j]$ protect all coefficients already fixed at stage $j$.
The supports of the blocks are pairwise disjoint, occur in successively
more distant symmetric bands, and all previous coefficients, including
zeros, remain fixed.

The estimates from Lemma~\ref{lem:phase-block} give
\begin{equation}\label{eq:phase-iteration-estimates}
 \norm{P_j}_1\le A+\sum_{k=1}^j\eta_k\le A+\eps,
 \qquad
 \HH(P_j)\le\theta\HH(P_{j-1})+\eta_j.
\end{equation}
Hence $\HH(P_j)\to0$, exactly as in~\eqref{eq:H-limit}. Stabilization
of the Fourier coefficients and the uniform $L^1$ bound produce a
weak-star limit $\nu\in M(\T)$ with
\begin{equation}\label{eq:phase-limit-transform}
 \wh\nu=A\one_{\{0\}}+\sum_{j\ge1}\wh Q_j.
\end{equation}
The gaps between the blocks give the same exact recovery identities
$V_{L_j}*\nu=P_j$ as in~\eqref{eq:circle-exact}. The proof of
Lemma~\ref{lem:ac-recovery} applies to complex measures, and
$\HH(P_j)\to0$ therefore implies $\nu\perp m_{\T}$. Since
$\wh\nu(0)=A$, the measure is nonzero and
$A\le\norm\nu\le A+\eps$.

Every nonzero coefficient in the $j$th block belongs to
$c_j\Omega_4$, proving~\eqref{eq:uniform-range}. The support estimate
\eqref{eq:phase-block-count}, together with
\eqref{eq:phase-c-choice}, gives, for every fixed $r>2$ and all large
$j$,
\[
 \sum_{n\in\supp\wh Q_j}|\wh Q_j(n)|^r
 \le\frac{A^2}{4}D_jc_j^{r-2}
 \le\frac{A^2}{4}D_jc_j^{1/j}
 \le\frac{A^2}{4}2^{-j}.
\]
Thus $\wh\nu\in\ell^r(\Z)$ for every $r>2$, and in particular
$\nu\in M_0(\T)$.

It remains to prove uniform convergence. Since $0<c_j<1$,
\eqref{eq:phase-c-choice} and~\eqref{eq:phase-block-U} yield
\begin{equation}\label{eq:phase-U-sum}
 \sum_{j\ge1}\norm{Q_j*Q_j}_U
 \le3A\sum_{j\ge1}D_jc_j
 \le3A\sum_{j\ge1}D_jc_j^{1/j}<\infty.
\end{equation}
Since $\norm h_\infty\le\norm h_U$ for every trigonometric
polynomial $h$, \eqref{eq:phase-U-sum} and the Weierstrass
$M$-test show that the series
\begin{equation}\label{eq:phase-square-series}
 f=A^2+\sum_{j\ge1}Q_j*Q_j
\end{equation}
converges uniformly to a continuous function. The spectral bands are
disjoint, so at every integer frequency at most one block is present and
\[
 \wh f(n)=\wh\nu(n)^2=\wh{\nu*\nu}(n).
\]

Let $j=j(N)$ be the first index for which
$\supp\wh Q_j$ is not entirely contained in $[-N,N]$. If $[-N,N]$
ends in the gap preceding the $j$th band, then
\[
 \norm{S_Nf-f}_\infty
 \le
 \sum_{k\ge j}\norm{Q_k*Q_k}_\infty
 \le
 \sum_{k\ge j}\norm{Q_k*Q_k}_U.
\]
If the cutoff meets the $j$th band, then all earlier bands are complete
and all later bands are absent, whence
\begin{equation}\label{eq:phase-partial-tail}
 \norm{S_Nf-f}_\infty
 \le
 2\norm{Q_j*Q_j}_U
 +\sum_{k>j}\norm{Q_k*Q_k}_U.
\end{equation}
As $N\to\infty$, one has $j(N)\to\infty$, and both bounds tend to zero
by~\eqref{eq:phase-U-sum}. Hence the Fourier series of $f$ converges
uniformly, and uniqueness of Fourier transforms gives
$\nu*\nu=fm_{\T}$.

The convergence cannot be absolute. Indeed, absolute convergence
would give
\[
 \sum_{n\in\Z}|\wh f(n)|
 =\sum_{n\in\Z}|\wh\nu(n)|^2<\infty,
\]
which would imply that $\nu$ has an $L^2$ density, contradicting the
fact that $\nu$ is nonzero and purely singular. Finally, the Rajchman
condition implies strong continuity by Lemma~\ref{lem:Rajchman-band}, and
$\nu*\nu\in L^1(\T)$ implies natural spectrum by
Lemma~\ref{lem:natural-spectrum}.
\end{proof}

\begin{remark}
The conclusion can be summarized as
\[
 \nu\perp m_{\T},
 \qquad
 \nu*\nu\in U(\T)\setminus A(\T),
\]
where $U(\T)$ and $A(\T)$ are the spaces defined above. Complex phases
are essential under the evenness requirement: no nonzero real even purely
singular measure has this property. The fourth roots of unity provide the
phases used in our construction. The theorem answers the circle question
for even complex measures when uniform convergence is understood in terms
of the ordinary symmetric Fourier partial sums
\cite[Remark~2.1, p.~16]{GuptaHare}.
\end{remark}

\subsection{Exact palettes versus Wiener--Pitt sets}

Theorems~\ref{thm:exact-torsion} and~\ref{thm:classification} show that
thinness of the value set alone does not prevent exact singular
realization. If $\Tor\wh G$ is infinite, every positive null palette,
with arbitrarily fast decay, is the exact range of a purely singular
Rajchman measure of mass $a_1$ and norm at most $a_1+\eps$.

This is compatible with \cite[Theorems~2 and~12]{OSW2026}. Fix a norm
bound $R>a_1$ and suppose that the palette satisfies the corresponding
decay assumptions in that work. Choosing $0<\eps<R-a_1$, the exact
realizations constructed here lie in the prescribed norm ball. The cited
theorems then imply that their convolution squares are absolutely
continuous and that the measures have natural spectra, while the measures
themselves remain purely singular. Thus the recent Wiener--Pitt theorem is
sharp with respect to the number of convolutions: the palette can force
regularity of $\mu*\mu$ and natural spectrum, but it cannot force
$\mu\in L^1(G)$.

On groups with finite $\Tor\wh G$, the counterpalette in
Theorem~\ref{thm:universal-counterpalette} gives the complementary
phenomenon. There are palettes for which no strongly continuous measure
can even attain all values, independently of singularity or natural
spectrum. The obstruction is not the rate of convergence alone; it is a
quantitative interaction between clusters of distinct values and the
number of large Fourier--Stieltjes coefficients allowed by the group.

\subsection{The open-range problem}

Ohrysko and Wojciechowski asked whether an open set of admissible
Fourier--Stieltjes values can force absolute continuity
\cite[Section~7, item~2]{OW2014}. Theorem~\ref{thm:open} answers this
negatively in every nonvacuous circle case, and more generally on the
groups stated there. The measure can have arbitrarily small norm, can be
Rajchman and strongly continuous, can have a natural spectrum, and can
have very regular second and third convolution powers, while retaining a
nonzero singular part. Hence none of those additional properties repairs
the failure of an open range condition to detect absolute continuity.

\subsection{Uniform convergence and the Gupta--Hare question}

Saeki \cite{Saeki1980} constructed a singular measure on $\T$ whose
convolution square has a uniformly convergent Fourier series. Gupta and
Hare \cite[Remark~2.1, p.~16]{GuptaHare} asked for a
reflection-invariant version and extended the construction to connected
compact groups. The displayed definition of the $u$-norm on p.~10 takes
the supremum over arbitrary finite subsets of the unitary dual, whereas
the proof on
pp.~12--13 uses initial segments of a fixed enumeration. These two notions
are different.

Proposition~\ref{prop:unconditional-obstruction} shows that the literal
finite-set interpretation is impossible for every nonzero singular measure
already on $\T$. A second obstruction applies to real even measures: if the
convolution square has a bounded density, then the measure has an $L^2$
density. Theorem~\ref{thm:uniform-square} therefore uses complex phases to
preserve evenness while turning coefficient squares into Rudin--Shapiro
signs. This gives uniform but nonabsolute convergence and shows that the
passage to complex measures is essential under the evenness requirement.

The result complements K\"orner's support and smoothness refinements
\cite{Korner2008}. The real-even flat-density assertion in
\cite{ElAbdalaoui2023} is not used here and is discussed in
Remark~\ref{rem:el-abdalaoui}, because under the reflection-invariant
interpretation it conflicts with Proposition~\ref{prop:real-even-obstruction}.

\subsection{A higher-convolution problem}

Section~7, item~3, p.~138 of \cite{OW2014} asks for value sets which
force absolute continuity at a prescribed convolution power for all
corresponding measures while still allowing examples whose earlier powers
remain singular. The construction of a single measure addresses a
different problem. In that direction, the basic blocks of the present
method give a singular $\mu$ with $\mu*\mu\in L^p$ for every
$1\le p<\infty$ and $\mu^{*3}$ in the Wiener algebra. The phase selector goes
further by arranging conditional uniform convergence of the square.

The local selector problem by itself is not an obstruction. Fix
$k\ge2$, put $\omega=e^{\pi i/k}$, and let complex numbers
$b_1,\ldots,b_{k-1}$ satisfy $\sum_{j=1}^{k-1}|b_j|\le1/2$. Define
\begin{equation}\label{eq:higher-selector-probabilities}
 p_r=\frac{1+2(-1)^r\RePart\left(
       \sum_{j=1}^{k-1}b_j\omega^{-rj}\right)}{2k},
 \qquad 0\le r<2k.
\end{equation}
These numbers are nonnegative, sum to one, and their sums over the even
and odd indices are both $1/2$. First define an ideal pair of random
variables $(F,\chi)$ by
\[
 \mathbb P\bigl(F=\omega^r,\ \chi=(-1)^r\bigr)=p_r,
 \qquad 0\le r<2k.
\]
Then $F^k=\chi$ and orthogonality on the cyclic group of order $2k$ gives
\[
 \E[\chi F^j]=b_j\qquad(1\le j<k).
\]
To realize a finite approximation on the Boolean cube, put
$\chi_{[d]}(s)=\prod_{j=1}^d s_j$. Within each parity class
$\{s:\chi_{[d]}(s)=\pm1\}$, round the target counts $2^dp_r$
separately, preserving the size $2^{d-1}$ of each class. Assigning the
corresponding values $\omega^r$ yields a selector on $\{-1,1\}^d$ such
that
\begin{equation}\label{eq:higher-selector-error}
 F^k=\chi_{[d]},
 \qquad
 \left|\wh{F^j}([d])-b_j\right|\le2k\,2^{-d}
 \quad(1\le j<k).
\end{equation}
Moreover, $|F|=1$, so Parseval and Cauchy--Schwarz automatically give
\[
 \sum_{S\subset[d]}|\wh F(S)|\le2^{d/2}.
\]
Thus neither finite-dimensional existence nor a finite Walsh
$\ell^1$ bound is the unresolved point.

What remains open for this method is an iterative scale problem. For a
prescribed $k\ge3$, one would need amplitudes $a_j$, numbers $M_j$ of nonzero block
coefficients, and selectors satisfying the moment conditions above while
simultaneously
providing: a summable $L^1$ error and a uniform Hellinger contraction;
a summable uniform-convergence or Wiener-algebra budget for the $k$th
convolution power; and exact recovery estimates certifying singularity of
each earlier power. A necessary scale condition for the latter is
\begin{equation}\label{eq:higher-necessary-scale}
 \sum_jM_ja_j^{2\ell}=\infty
 \qquad(1\le\ell<k),
\end{equation}
for blocks having $M_j$ nonzero coefficients of modulus $a_j$, because
$\sum_jM_ja_j^{2\ell}<\infty$ would put the Fourier transform of the
$\ell$th convolution power in $\ell^2$ and hence make that power
absolutely continuous. The present block estimate has the scale
$M_j\lesssim D_ja_j^{-2}$, and the iteration chooses
$D_ja_j^{1/j}$ summably small. Consequently
$\wh\mu\in\bigcap_{r>2}\ell^r$, so $\mu*\mu$ is already absolutely
continuous, independently of the phases. The open issue is to satisfy the displayed iterative
requirements simultaneously, not to construct the local selectors alone.

\subsection{Further directions}

The classification leaves several natural questions.
\begin{enumerate}[label=\textup{(\arabic*)}]
\item Characterize the individual palettes exactly realizable by an
atomless singular measure on $\T$. The cost condition
$\Lambda_K(p)\to0$ is sufficient but is not known to be necessary.
\item Determine whether the norm $a_1+\eps$ in exact Rajchman
realization on compact abelian groups whose dual has an infinite torsion
subgroup can be replaced by the optimal value $a_1$ without losing pure
singularity.
\item Classify exact realization by positive singular measures. The
signed classification in Section~\ref{sec:exact-palettes} does not
imply a positive one.
\item Extend the palette classification to complex compact palettes
with zero as their only accumulation point. The ordered positive
structure is used essentially in both telescoping and bounded-cost
tranche forcing.
\end{enumerate}

\section*{Use of artificial intelligence}

The author received access to OpenAI tools through the ChatGPT for Academic
Researchers program. AI-assisted tools were used during the preparation
of this manuscript for mathematical verification, literature searches,
and editorial support. The author has independently verified the content
and takes full responsibility for its accuracy and conclusions.

\end{document}